\documentclass[11pt]{article}

\usepackage{amsfonts,amssymb, latexsym, amsmath, amsthm,mathrsfs, verbatim,geometry}
\usepackage{graphics}
\usepackage{slashed}
\usepackage{url,color}
\usepackage{enumerate}
\usepackage{esint}
\usepackage{pifont}
\usepackage{upgreek}
\usepackage{times}
\usepackage{calligra}
\usepackage{graphicx}
\usepackage{caption}

\usepackage{tikz}
\usetikzlibrary{calc}
\usetikzlibrary{calc}
\usetikzlibrary{arrows.meta}
\usetikzlibrary{hobby}
\usepackage{tkz-euclide}

\usepackage{caption}

\usepackage{hyperref}
\hypersetup{
    colorlinks,
    citecolor=red,
    filecolor=green,
    linkcolor=blue,
    urlcolor=black
}

\def\R{D}
\def\C{\mathbb C}

\newcommand\J{\mathscr J}
\def\A{\mathscr{A}}

\def\D{\mathcal{D}}
\def\be{\begin{equation}}
\def\ee{\end{equation}}
\def\bea{\begin{eqnarray}}
\def\eea{\end{eqnarray}}
\def\beas{\begin{eqnarray*}}
\def\eeas{\end{eqnarray*}}

\def\l{\lambda}
\def\pa{\partial }

\def\r{ r\partial_r}

\def\l{\lambda}

\def\lv{\left\vert}
\def\rv{\right\vert}

\def\w{{\bf w}}
\def\d{{\bf d}}

\def\chid{\chi_{\text{dust}}}

\def\bcr{\begin{color}{red}}
\def\bcb{\begin{color}{blue}}
\def\bcg{\begin{color}{green}}
\def\bcv{\begin{color}{violet}}
\def\ec{\end{color}}

\def\k{\varepsilon}
\def\kl{\epsilon}

\def\F{J}

\def\xmin{x_{\text{min}}}
\def\xmax{x_{\text{max}}}

\def\wl{{\bf w}}

\def\RY{d}
\def\vY{w}

\def\yms{Y^{\text{ms}}}

\def\MS{\mathcal{M}\mathcal{S}_\k}

\def\u{p}
\def\Om{\Omega}

\def\uC{\underline{\mathcal C}}
\def\C{\mathcal C}

\def\ga{\gamma}
\def\om{\omega}

\newcommand{\vertiii}[1]{{\left\vert\kern-0.25ex\left\vert\kern-0.25ex\left\vert #1 
    \right\vert\kern-0.25ex\right\vert\kern-0.25ex\right\vert}}

\newtheorem{theorem}{Theorem}[section]

\newtheorem{lemma}[theorem]{Lemma}
\newtheorem{remark}[theorem]{Remark}

\title{Star dynamics: collapse vs. expansion}

\author{Mahir Had\v zi\'c\thanks{Department of Mathematics, University College London, London UK, Email: m.hadzic@ucl.ac.uk.}}

\date{}

\begin{document}

\maketitle

{\centering\footnotesize Dedicated to Constantine M. Dafermos on the occasion of his 80-th birthday.\par}

\abstract{
We review a series of recent results on global dynamic properties of radially symmetric self-gravitating compressible Euler flows, which naturally arise in the mathematical description of stars. We focus on the role of scaling invariances and how they interact with  nonlinearities to generate imploding finite-time singularities as well as expanding star solutions, arising from smooth initial data. This review paper is based on joint works with Y. Guo, J. Jang, and M. Schrecker.}

\tableofcontents

\section{Introduction}

First systematic mathematical studies of the global properties of self-gravitating fluids  focused on steady states (i.e. time-independent solutions), which have played a central role in both astrophysics and mathematics literature, see for example the classical textbook of Cahndrasekhar~\cite{Chandrasekhar1939}. It is in particular natural to use stationarity to model objects such as the neutron stars and white dwarves, see the discussion in~\cite{LiYa1987}. Different thermodynamic properties of the gas in the star are then encoded in a choice of the so-called equation of state, which can be rather complicated in realistic situations - we refer the interested reader to~\cite{Chandrasekhar1939,HaThWaWh, Th1966, Chandrasekhar1984, ZeNo, Straumann2013} for more astrophysical background on the topic. Mathematically, the study of stationary stars and their stability properties is a beautiful topic with many unsolved open questions.

Our goal here is to move away from stationarity (and the associated stability analysis) and study time-dependent regimes wherein basic features of the physically important phenomena such as stellar collapse or supernova explosions can be captured in a mathematically rigorous way. We are in part motivated by the astrophysical expectation that the end-state of the gravitational collapse of a star should take on a self-similar form to leading order. This informal statement is sometimes referred to as the self-similarity hypothesis~\cite{CaCo1999,Ha1998,Ha2003,Ha2004, HaMa2001}. It naturally leads to the basic question of the existence of physically meaningful self-similar solutions representing for example collapse - this topic is at the centre of this review.

Many of the mathematical difficulties associated with self-gravitating fluids originate in a very nonlinear and nonlocal interaction between the pressure and the gravity. The former typically encourages dispersion, while the latter is clearly an attractive, focusing force. This antagonism, among other things, allows for the existence of steady states mentioned above -- this is the case when these two forces are in perfect balance. However, it also leads to a number of different regimes, including most notably collapse and expansion.

The basic model of a self-gravitating dynamically evolving star is given by the compressible Euler-Poisson system, from now on the EP-system. The system takes the form
\begin{align}
\pa_t\rho+\text{div}(\rho{\bf u}) & = 0 \label{E:EPCONT}\\
\rho D_t {\bf u} + \nabla p + \rho \nabla \phi & = 0 \\
\Delta\phi & = 4\pi\rho,\label{E:EPPOISSON}
\end{align}
where the unknowns are the density $\rho$, the pressure $p$, the gravitational potential $\phi$, and the velocity 3-vector ${\bf u}$. To ensure that we model an isolated body, we demand asymptotic flatness encoded in the assumption $\lim_{|x|\to\infty}\phi(x)=0$. As written,~\eqref{E:EPCONT}--\eqref{E:EPPOISSON} is underdetermined, with pressure $p$ an unknown one too many.  A standard way of closing the system is to impose an equation of state, which relates the pressure to the density. In this article, we are interested in the class of polytropic equations of state
\begin{align}\label{E:EOS}
p = \k \rho^\gamma,  \ \ 1\le \gamma< 2, \ \k>0.
\end{align}
The polytropic index $\gamma$ parametrises the above family of problems and serves as a natural criticality parameter, see Section~\ref{S:SCALING}.
The case $\gamma=1$ is sometimes referred to as the  isothermal case and plays a particularly important role in the description of stellar collapse.

The EP-system has received a lot of attention in the astrophysics literature. A classical topic in this respect is the existence and (mostly linear) stability 
analysis of a well-known class of radial steady state solutions known as Lane-Emden stars~\cite{Chandrasekhar1939,ZeNo,BiTr}, with many contributions in both mathematics and physics literature, see for example
the introduction to~\cite{HaJa2016-1} for the literature overview. At a nonlinear level, rigorous description of the phase space around the Lane-Emden stars remains an outstanding open problem. This includes the questions of the asymptotic-in-time behaviour, the existence of (quasi-) periodic solutions near the Lane-Emden stationary profiles, and the formation of shocks from smooth and small perturbations. There are only a few rigorous results available: conditional nonlinear stability by Rein~\cite{Rein2003} (see also the work of Luo and Smoller~\cite{LuSm}) and rigorous nonlinear instability results in the range $\frac65\le \ga<\frac43$ by Jang~\cite{Jang2008,Jang2014}.\footnote{Although viscous effects are not considered in this review, we mention that the stability properties of Lane-Emden stars are better understood in the presence of viscosity, see~\cite{JaTi2013,LXZ1,LXZ2,HoLuZh2018}.}

In this review we focus on the nonlinear dynamics in a region of the phase space that is far from stationarity. Section~\ref{S:NEWTON} is devoted to the Newtonian theory, mostly the EP-system. We study two interesting dynamic scenarios which have a basis in astrophysical observations -  stellar collapse and stellar expansion. The former is mathematically captured through finite-time singularity formation driven by the density blow up, also referred to as implosion. Such singularities are therefore not shocks, as shocks are characterised by the formation of discontinuities in the field variables and give the archetype singular behaviour associated with compressible flows~\cite{Da2010}.  

The common thread to all of our results is the decisive role of scaling invariances. In Sections~\ref{S:LAGR} and~\ref{S:SCALING} we present the Lagrangian formulation of the problem and explain the scaling symmetries respectively. We then show in Section~\ref{SS:EXEX} that one can find examples of global-in-time solutions to compressible Euler and Euler-Poisson flows which are driven by the expansion of the fluid (star) support - these are the main results of~\cite{HaJa2016-1} and~\cite{HaJa2016-2}. Working in the opposite direction, 
we also prove that there exist at least two types of collapsing stars in the Newtonian case. In Section~\ref{SS:DUST} we explain the construction of collapsing stars from~\cite{GHJ2021a} that on approach to singularity behave like dust collapsing stars to the leading order. In contrast to such dust-like behaviour, in Section~\ref{SS:NEWTONIANCOLLAPSE} we explain our recent construction~\cite{GHJ2021b,GHJS2021} of exactly self-similar collapsing flows, discovered numerically by Larson and Penston in 1969~\cite{Larson1969,Penston1969} when $\ga=1$, and by Yahil in 1983~\cite{Yahil83} for certain ranges of $\ga>1$. The central mathematical challenge in our analysis is the presence of so-called sonic points. They encode natural singularities for the underlying dynamical system and correspond geometrically to the boundary of the backward acoustical cone emanating from the singularity. In a very recent work Merle, Rapha\"el, Rodnianski, and Szeftel~\cite{MRRS1} proved the existence of $C^\infty$ imploding self-similar solutions for the Euler system (without gravity) featuring a sonic point.

In the context of general relativity, one must couple the Einstein field equations to the Euler equations. In this setting the gas density, pressure, and the 4-velocity $u^\mu$ evolve simultaneously with the underlying Lorentzian manifold $(\mathcal M,g)$ through the Einstein equations:
\begin{align}\label{E:EINSTEIN}
R_{\mu\nu}-\frac12 R g_{\mu\nu} = 8\pi T_{\mu\nu}, \ \ \mu,\nu=0,1,2,3.
\end{align}
Here, $T_{\mu\nu}$ denotes the energy-momentum tensor which takes the form
\begin{align}\label{E:EMTENSOR}
T_{\mu\nu} = (\rho+p)u_{\mu}u_\nu + p g_{\mu\nu},
\end{align}
and the $4$-velocity is a future directed timelike vector satisfying the normalisation condition
\begin{align}\label{E:NORMALISATION}
g_{\mu\nu}u^\mu u^\nu = -1.
\end{align}
We refer to equations~\eqref{E:EINSTEIN}--\eqref{E:NORMALISATION} as the Einstein-Euler system, from now on the EE-system.
Here we shall assume the equation of state~\eqref{E:EOS} with $\gamma=1$. In this case the constant $\k$ corresponds to the square of the speed of sound
and to respect causality we assume $0\le \k\le1$. 

Under the right rescaling, the EP-system can be thought of as the Newtonian limit of the EE-system. The same questions as in the case of the EP-system apply to the EE-system, but due to the structure of the Einstein field equations and its interaction with the 
compressible Euler, much less is known. The analogues of the Lane-Emden stars have been known to exist since 1930's, going back to the classical works of Tolman, Volkov, and Oppenheimer~\cite{To1934-2,OpVo1939}. Here the equation of state takes the form~\eqref{E:EOS} with $\gamma=1$ at high densities, and with $\gamma\ge1$ at small densities close to the vacuum boundary. Formal linear stability analysis was initiated by Chandrasekhar in 1960's~\cite{Chandrasekhar1964} and the related turning point principle was formulated by Harrison, Thorne, Wakano, Wheeler, Zeldovitch et al.~\cite{HaThWaWh,Th1966,Ze1963}. A rigorous linear (in)stability theory and the proof of the turning point principle were given recently by the author jointly with Z. Lin and G. Rein~\cite{HaLiRe2021,HaLi2021}. 

In the absence of pressure, the relativistic fluid is often referred to as dust. In early 1930s Tolman~\cite{To1934} and Lema\^itre~\cite{Lemaitre1933} solved the Einstein-dust system in radial symmetry (see also the later work of Bondi~\cite{Bondi1947}), thus obtaining an infinite-dimensional family of solutions to the Einstein-dust system. Such solutions contain open families of initial data that lead to finite-time gravitational collapse, at first instance measured by the blow-up of the dust mass-energy density. In their seminal work from 1939, Oppenheimer and Snyder~\cite{OpSn1939} focused on a subclass of Tolman solutions characterised by the space-homogeneity of the dust density. By analysing the causal structure of the ensuing spacetime they laid the foundation for the
modern notion of a black hole. However, in 1984 Christodoulou~\cite{Ch1984} examined the causal structure of the remaining, space-inhomogeneous (and therefore generic) Tolman spacetimes and came to a surprising conclusion. Instead of black holes, such spacetimes contained so-called naked singularities, which loosely speaking, emit signals that can be received by asymptotic observers. Such singularities are conjectured to be a non-generic feature of the theory of general relativity, and this is often mathematically codified via cosmic censorship hypotheses - we refer the reader to~\cite{Ch1999b, DaLu2017, RoSR2019} for a rigorous discussion of this topic. In the context of radially symmetric matter models, Christodoulou proved the existence of a class of naked singularity spacetimes for the Einstein--scalar field model~\cite{Ch1994}. These solutions are self-similar in a suitable sense. In line with the weak cosmic censorship conjecture, he also proved that such naked singularities are unstable~\cite{Ch1999b}. Much more recently, Rodnianski and Shlapentokh-Rothman~\cite{RoSR2019,SR2022} proved the existence of naked singularity spacetimes solving the Einstein-vacuum system of equations.

In view of the above discussion, the conclusions of~\cite{Ch1984} can be taken as a motivation for study of the causal structure of implosion for more realistic fluid models which include pressure. In Section~\ref{S:NAKED} we finally present our recent result~\cite{GHJ2021c} on the formation of self-similar imploding singularities
for the EE-system with $\ga=1$ and $0<\k\ll1$ in~\eqref{E:EOS}, and thereby make rigorous the pioneering predictions of Ori and Piran~\cite{OP1987,OP1988,OP1990}. The resulting spacetimes are referred to as the relativistic Larson-Penston solutions. Again, here we are dynamically far from stationarity and the constructed spacetime, as we will show, contains a naked singularity. The latter statement requires an in-depth description of the causal structure of the spacetime and this interaction between the Lorentzian- and the acoustical geometry associated with the propagation of sound waves is one of the most interesting features of this problem. 

\medskip 

{\bf Acknowledgments.}
The author acknowledges the support of the EPSRC Early Career Fellowship EP/S02218X/1.
He thanks Y. Guo, J. Jang, and M. Schrecker for fruitful discussions.

\section{Self-gravitating Newtonian gases}\label{S:NEWTON}


\subsection{Lagrangian point of view} \label{S:LAGR}

In the absence of vacuum regions, local well-posedness for the EP-system can be shown using for example the classical symmetrisation techniques~\cite{Da2010} for hyperbolic conservation laws combined with an energy method. In the presence of vacuum the situation is more complicated, see for example~\cite{MaUkKa,MaPe1990}. For more references on the topic and global-in-time weak solution theories in radial symmetry, we refer the reader to the recent work~\cite{ChHeWaYu}.

In the context of classical compressible fluid dynamics it is often advantageous to use 
the Lagrangian (or comoving) coordinates, which turn out to be particularly well-suited to the presence of vacuum regions. If we denote by $\Omega(t)\subset \mathbb R^3$ the support of the gas, the basic unknown is the flow map 
$\eta(t,\cdot)$ defined on some reference domain $\mathfrak M\subset \mathbb R^3$
so that
$\eta(t,\cdot):\mathfrak M\to \Omega(t)$ solves
\begin{align}
\pa_t\eta(t,x) & = {\bf u}(t,\eta(t,x)), \ \ x\in \mathfrak M, \label{E:ETADEF}\\
\eta(0,x) & = \eta_0, \ \ x\in\mathfrak M.\label{E:ETADEFINITIAL}
\end{align}
Here $\eta_0:\mathfrak M\to \Omega(0)$ is the labelling map and its choice fixes the labelling gauge freedom in the problem; recall that ${\bf u}(t,\cdot):\Omega(t)\to\mathbb R^3$ is the velocity vectorfield associated with the gas flow. In the case of a free-moving surface $\pa\Omega(t)$, Lagrangian coordinates can be used to pull the problem back to the fixed reference domain $\mathfrak M$. To do so, we introduce the new unknowns
\begin{align}
{\bf v} : = {\bf u}\circ \eta \ \ &\text{ (Lagrangian velocity)},\label{E:VDEF}\\
\A : = [D\eta]^{-1} \ \ &\text{ (Inverse of the Jacobian matrix)},\\
\J  : = \det [D\eta] \ \ &\text{ (Jacobian determinant)}.
\end{align}
For notational brevity, it is useful to introduce
\begin{align}
\alpha:=\frac1{\ga-1}.
\end{align}
A standard calculation~\cite{JaMa2015} shows that the the compressible Euler equations take the form
\begin{align}\label{E:EULERL}
w^\alpha \pa_{tt}\eta_i + \left(w^{1+\alpha}\A^k_i\J^{1-\gamma}\right),_k = 0, \ \ i=1,2,3, \ \ x\in \mathfrak M,
\end{align}
where 
\begin{align}\label{E:WDEF}
w(x) :=  \left(\rho_0(\eta_0(x)) \J_0(x)\right)^{\gamma-1}
\end{align}
corresponds to the initial gas enthalpy. Here we have assumed $\k=1$ in~\eqref{E:EOS} as the exact value does not play a role for the considerations in this section. We note that $w(x)=\rho_0(x)^{\ga-1}$ if $\mathfrak M = \Omega(0)$ and $\eta_0=\text{Id}$.
In the presence of a self-induced gravitational field, the Lagrangian description of the Euler-Poisson system~\eqref{E:EPCONT}--\eqref{E:EPPOISSON}
takes the form
\begin{align}
w^\alpha \pa_{tt}\eta_i + \left(w^{1+\alpha}\A^k_i\J^{1-\gamma}\right),_k + \A^k_i\psi,_k &= 0, \ \ i=1,2,3, \ \ x\in \mathfrak M, \label{E:EPL}\\
\A^k_i(\A^j_i\psi,_j),_k &= 4\pi w^\alpha \J^{-1}, \label{E:EPPOISSONL}
\end{align}
where $\psi=\phi\circ\eta$ is the Lagrangian pull-back of the gravitational potential $\phi$.
As we can see from~\eqref{E:EULERL} and~\eqref{E:EPCONT}--\eqref{E:EPPOISSON}, both systems of equations are quasilinear second order wave-like systems.
The precise nature of the flow and the role of the enthalpy in the context of well-posedness for the free surface flows will become apparent in later sections.

 If the flow is radially symmetric, we can write the flow map in the form
\[
\eta(t,x) = \chi(t,r) x,
\]
where $r=|x|$ and $\chi(t,\cdot)$ is a non-negative function. It is straightforward to check~\cite{Jang2014} that~\eqref{E:EULERL} takes the form
\begin{align}\label{E:EULERRADIALLAGR}
\chi_{tt} + \frac{\chi^2}{w^\alpha r}\pa_r\left(w^{1+\alpha}\J[\chi]^{-\gamma}\right) = 0,
\end{align}
while the radial Euler-Poisson system~\eqref{E:EPL}--\eqref{E:EPPOISSONL} takes the form
\begin{align}\label{E:EPRADIALLAGR}
\chi_{tt} + \frac{\chi^2}{w^\alpha r}\pa_r\left(w^{1+\alpha}\J[\chi]^{-\gamma}\right)  + \frac{G(r)}{\chi^2} = 0,
\end{align}
where $G(r)$ is the average density expressed in Lagrangian coordinates and $\J[\chi]$ the Lagrangian determinant
\begin{align}\label{E:AVERAGEDENSITY}
G(r) &= \frac1{r^3}\int_0^r 4\pi \rho_0(\chi_0(s)) \J_0(s) s^2\,ds,  \\
\J[\chi] & = \chi^2(\chi+r\pa_r\chi). \label{E:DETFORMULA}
\end{align}

\subsection{Scaling invariance}\label{S:SCALING}
It is a simple algebraic exercise to verify that for any $\l>0$ the transformation
\begin{align}\label{E:SCALING}
\rho \mapsto \lambda^{-\frac{2}{2-\gamma}} \rho(s,y),  \ \
\mathbf{u}(t,x) \mapsto \lambda^{-\frac{\gamma-1}{2-\gamma}} \mathbf{u}(s,y),
\ \ 
\phi(t,x) \mapsto \lambda^{-\frac{2\gamma-2}{2-\gamma}}\phi(s,y)
\end{align}
keeps the EP-system invariant, 
where the time and space scale according to
\begin{align}\label{E:SCALING2}
s = \frac{t}{\lambda^{\frac{1}{2-\gamma}}}, \ \ y = \frac{x}{\lambda}.
\end{align}
The polytropic index $\gamma$ therefore acts as a natural criticality parameter, and it is readily verified that the 
total mass
\be\label{E:MASSDEF}
M[\rho] : = \int_{\mathbb R^3} \rho \,dx,
\ee
and the total energy
\be\label{E:ENERGYDEF}
E[\rho,{\bf u}] := \int_{\mathbb R^3} \left(\frac12 \rho|{\bf u}|^2 + \frac1{\gamma} \rho^\gamma - \frac1{8\pi}|\nabla\phi|^2\right)\,dx,
\ee
are formally preserved exactly when $\gamma=\frac43$ and $\gamma=\frac65$ respectively.


\subsection{Existence via expansion}\label{SS:EXEX}


It is well-known that the tendency of fluid particles to expand and move away from each other can generate dispersion,
which counteracts possible focussing of the characteristics. This mechanism therefore has the potential to generate solutions that exist globally-in-time.
As it is expected for gravity to be a subdominant effect in such a situation, it is instructive to look for global-in-time
solutions to the compressible Euler system
\begin{align}
\pa_t\rho+\text{div}(\rho{\bf u}) & = 0 \label{E:ECONT}\\
\rho \left(\pa_t+{\bf u}\cdot\nabla\right){\bf u} + \nabla p  & = 0  \label{E:EMOM},
\end{align}
where we assume the equation of state~\eqref{E:EOS} with $\k=1$. 

First solutions of this flavour were constructed in the important works of Ovsyannikov~\cite{Ov1956} and Dyson~\cite{Dyson1968}. They belong to a class of so-called affine motions which will be discussed below. Expansion is also a mechanism for global existence in the works of Grassin~\cite{Grassin98}, Serre~\cite{Se1997}, and Rozanova~\cite{Ro}. 
In a remarkable work~\cite{Sideris2017}, Sideris exhibited an explicit class of expanding solutions to the Euler system~\eqref{E:ECONT}--\eqref{E:EMOM} that exist for all
$t>0$ and retain compactly supported density.
The Sideris flows are most easily described using the comoving description of the flow. He
restricts his attention to a special type of solutions, known as affine flows, which by definition take the form
\begin{align}\label{E:AFFINE}
\eta(t,x) = A(t) x, \ \ A\in \text{GL}^+(3,\mathbb R),
\end{align}
where we recall~\eqref{E:ETADEF}.
The compressible Euler system then reduces to 
\begin{align}\label{E:AFFINEB}
\ddot A(t) x  + \frac\gamma{\gamma-1}A(t)^{-\top} \det(A(t))^{1-\gamma} \nabla_x w & = 0, \ \ \gamma>1,
\end{align}
where $A^{-\top}$ stands for the transpose of the inverse of the matrix $A$ and $w = \rho_0^{\gamma-1}$ is the initial enthalpy.
The idea is to separate variables. Given a $\delta\in\mathbb R$ and 
$(A_0,A_1)\in\text{GL}^+(3,\mathbb R)\times\mathbb M^3$ we seek to solve the system
\begin{align}
\ddot A(t) & = \delta \det(A(t))^{1-\gamma} A(t)^{-\top}, \label{E:BEQN} \\
A(0) & = A_0, \ \ \dot A(0) = A_1, \label{E:BINITIAL} \\
\nabla_x w  & = - \delta \frac{\gamma-1}{\gamma} x.\label{E:FZEROEQN}
\end{align}
It is clear that any triple $(A, w,\delta)$ solving~\eqref{E:BEQN}--\eqref{E:FZEROEQN} gives a solution of~\eqref{E:AFFINEB}.
Equation~\eqref{E:FZEROEQN} can be solved explicitly and
gives 
\begin{align}\label{E:wform}
w(x) = \frac{\delta (\gamma-1)}{2\gamma}\left(1 - |x|^2\right), \ \ x\in B_1(0),
\end{align}
whereas, by~\cite{Sideris2017}, for 
any $\delta>0$ and $(A_0,A_1)\in\text{GL}^+(3,\mathbb R)\times\mathbb M^3$ 
there exists a unique, smooth, global-in-time solution to the initial value problem~\eqref{E:BEQN}--\eqref{E:BINITIAL}.

\begin{figure}
\begin{center}
\begin{tikzpicture}[>=stealth]
\def\wx{3}      
\def\wy{1.3}    
\def\wz{0.8}   

\draw (0,0) ellipse (\wx cm and \wy cm);

\begin{scope}
  \clip (-\wx-0.1, 0) rectangle (\wx+0.1, \wy);
  \draw[dashed] (0, 0) ellipse (\wx cm and \wz cm);
\end{scope}

\begin{scope}
  \clip (-\wx-0.1, 0) rectangle (\wx+0.1, -\wy);
  \draw (0, 0) ellipse (\wx cm and \wz cm);
\end{scope}

\def\r{1.2}    
\def\s{1.5}    
\foreach \x in {90, 180, ..., 360} {
  \draw[-{Latex[width=2mm]}] ({sin(\x)*\wx*\r}, {cos(\x)*\wy*\r}) -- ({sin(\x)*\wx*\s}, {cos(\x)*\wy*\s});
}
\end{tikzpicture}
\end{center}
\caption{Sideris expanding ellipsoids}
\label{F:SIDERIS}
\end{figure}
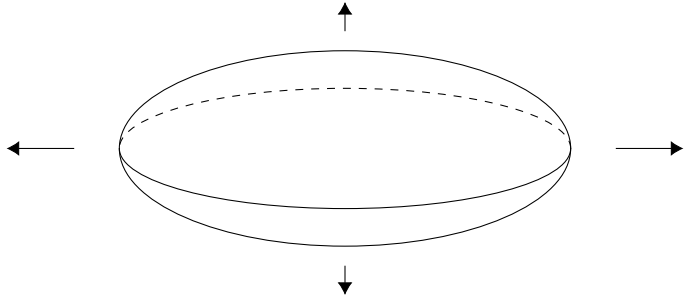

We may similarly look for the affine motions solving EP$_\gamma$ system. This leads us to a famous class of solutions discovered in 1980 by Goldreich and Weber~\cite{GoWe1980}, see~\cite{Makino92,FuLin,DengXiangYang} for subsequent mathematical work. Due to constraints imposed by the presence of gravity, the space of such solutions is considerably more rigid. For scaling reasons, we must assume that $\gamma=\frac43$ and we further work in radial symmetry, i.e. we assume that $A(t) = \lambda(t)\text{Id}$. The resulting ODE for $\l(t)$ 
reads
\begin{align}\label{E:GWODE}
\ddot\l(t)\l(t)^2 & = \delta, \ \ \delta\in\mathbb R, \\
\l(0)=\l_0, \ \ \dot\l(0)& =\l_1 \label{BClambda}
\end{align}
and the associated enthalpy satisfies the following generalised Lane-Emden equation
\begin{align}\label{E:GLE2}
\Delta_r w(r) + \pi w(r)^3 =- \frac34 \delta.
\end{align}

The Goldreich-Weber solutions split into two groups: the expanding star solutions (not unlike the Sideris flows above) and the collapsing star solutions, which are
very different from the Sideris flows. There also exists a purely Eulerian approach to the derivation of these affine flows developed by the author and Jang~\cite{HaJa2016-2}. While less direct, it clarifies the precise role of symmetries of the Euler equation responsible for the existence of such special structures and contextualises them within the broader setting
of dispersive PDE.

Both classes of affine motions described above feature a vacuum boundary separating regions of positive and vanishing pressure. More importantly, in both
cases the fluid pressure is not a smooth function at the vacuum boundary. This lack of smoothness is characterised by the behaviour of the fluid enthalpy~\eqref{E:WDEF} 
at the vacuum boundary, which both in the case of the Sideris- and the Goldreich-Weber flows reads
\be\label{E:PHYSICALVACUUM}
w(x)\asymp_{x\to\pa\Omega} \ \text{dist}(x,\pa\Omega).
\ee

Condition~\eqref{E:PHYSICALVACUUM} is known as the physical vacuum condition. In two independent works by Jang and Masmoudi~\cite{JaMa2015}
and Coutand-Shkoller~\cite{CoSh2012} it was shown that the free surface Euler equations are well-posed locally-in-time precisely under the assumption of physical vacuum on the initial data. These works give the natural mathematical framework for the study of dynamic properties of solutions to Euler flows. The techniques in~\cite{JaMa2015, CoSh2012} can be adapted to yield well-posedness for the Euler-Poisson system in the presence of physical vacuum~\cite{GuLe2016,HaJa2017}. For an Eulerian approach to the well-posedness in the presence of physical vacuum see the recent work~\cite{IfTa2020}.

\subsubsection{Nonlinear stability}

The natural question is whether the Sideris and the expanding Goldreich-Weber solutions are dynamically stable in some 
reasonable topology. We answered this question to the affirmative, taking advantage of the expansion of the underlying solution.

\begin{theorem}[\cite{HaJa2016-1}, informal statement]\label{T:GWGLOBAL}
The expanding Goldreich-Weber flows are nonlinearly stable in the class of radially symmetric perturbations.
\end{theorem}

\begin{theorem}[\cite{HaJa2016-2}, informal statement]\label{T:SIDERISGLOBAL}
Sideris flows are nonlinearly stable for $1<\ga\le\frac53$.
\end{theorem}

The main idea of the proofs of Theorems~\ref{T:SIDERISGLOBAL} and~\ref{T:GWGLOBAL} is to capture the stabilising effect of the background expansion in the Lagrangian framework.
If $(\rho,{\bf u})$ is an Eulerian description of a flow which starts off ``close" to a fixed Sideris motion $(t,x)\mapsto A(t)x$, we consider the associated flow map $\zeta$, i.e. $\pa_t\zeta = {\bf u}(t,\cdot)\circ\zeta$. The idea is to renormalise the map $\zeta$ by the background flow $A(t)x$ in a suitable way - we introduce the vectorfield 
$\uptheta$ through the relation
\begin{align}
A^{-1}(t) \zeta(t,x) = \text{Id} + \uptheta(t,x).
\end{align}
To give a succinct explanation of the key stabilisation mechanism, let us choose a particularly convenient subclass of diagonal Sideris flows of the form
\begin{align}
A(t) = \text{diag} (\l(t),\l(t),\l(t)),
\end{align}
where according to~\eqref{E:BEQN}--\eqref{E:BINITIAL} $\l(t)$ satisfies the second order ODE
\be\label{E:AFFINESIMPLE}
\ddot\l(t)\l(t)^{3\ga-2} = \delta,\ \ \l(0)=1, \ \ \dot\l(0)=\l_1\in\mathbb R.
\ee
A routine exercise in ODE shows that for any $\delta>0$ there exists a global solution to~\eqref{E:AFFINESIMPLE} and $\l(t)\asymp_{t\to\infty}t$. In other words, the associated
affine motion expands at a linear rate in $t$. To highlight the stabilising effect of this expansion we switch to a logarithmic time-scale by introducing a new time variable $\tau$ through
$\frac{d\tau}{dt}=\frac1{\l(t)}$. It is another routine exercise to derive the PDE satisfied by the perturbation $\theta$:
\begin{align}\label{E:THETAPDE}
 \l(\tau)^{3\ga-3}\left(\pa_{\tau\tau}\theta_i +  \frac{\pa_\tau\l}{\l}\pa_\tau\theta_i\right) +\theta_i + \frac{\gamma}{\gamma-1} w^{-\alpha} \pa_k\left( \A^k_i \J^{-\frac{1}{\alpha}} - \delta^k_i\right) =0, \ \ i=1,2,3.
\end{align}
In self-similar variables the linear-in-$t$ growth of $\l(t)$ translates into the asymptotic property $\l(\tau)\asymp_{\tau\to\infty} e^\tau$. In particular $\frac{\pa_\tau\l}{\l}\approx_{\tau\to\infty} 1$. It is therefore illustrative to consider the simpler model equation
\begin{align}
 &e^{(3\ga-3)\tau}\left(\pa_{\tau\tau}\theta_i +  \pa_\tau\theta_i\right) +\theta_i + F_i[\theta] =0, \ \ i=1,2,3,\label{E:THETAPDE2} \\ 
 & F_i[\theta] = \frac{\gamma}{\gamma-1} w^{-\alpha} \pa_k\left( w^{1+\alpha} \left(\A^k_i \J^{-\frac{1}{\alpha}} - \delta^k_i\right)\right), \ \ i=1,2,3. \label{E:ELLIPTIC}
\end{align}
Upon evaluating the inner product of~\eqref{E:THETAPDE2} with $\pa_\tau\theta^i$ in the weighted $L^2_\alpha(\Omega):=L^2(\Omega, w^\alpha\,dx)$ space,
the first three terms on the left-hand side of~\eqref{E:THETAPDE2} give the contribution
\begin{align}
&\left( e^{(3\ga-3)\tau} \left(\pa_{\tau\tau}\theta_i +  \pa_\tau\theta_i\right) +\theta_i\,,\pa_\tau\theta^i\,\right)_{L^2_\alpha} \notag\\
&= \frac12\pa_\tau\left(e^{(3\ga-3)\tau}\|\pa_\tau\theta\|_{L^2_\alpha}^2 + \|\theta\|_{L^2_\alpha}^2\right) + \frac{5-3\ga}{2}e^{(3\ga-3)\tau}\|\pa_\tau\theta\|_{L^2_\alpha}^2 \label{E:DAMPING}
\end{align}
Equation~\eqref{E:DAMPING} gives a clear stabilising effect in the regime $\ga\le\frac53$ and we refer to it as ``damping", although it should be more precisely be referred to as dispersion. Since $\tau$-derivatives commute with the spatial ones, it is clear that the analogue of the damping term above will appear in a high-order energy framework based on commuting~\eqref{E:THETAPDE2} with spatial derivatives. The nonlinear term $F_i[\theta]$~\eqref{E:ELLIPTIC} is on the other hand well-known in the context of free-boundary vacuum problems. If $\partial$ denotes the highest-order spatial derivative that we commute $F_i[\theta]$ with, upon multiplying the commuted equation by $\pa \theta^i$, our goal is to extract a (hopefully) positive-definite coercive energy-like term upon integrating-by-parts. The story is however not so simple since the presence of the weight $w$ degenerates the operator on approach to the boundary $\pa\Omega$. This issue is at the heart of the complications associated with the local-in-time well-posedness. For our purposes, the resolution came about through a deep idea of Jang and Masmoudi~\cite{JaMa2015} where they implemented a particular weighting scheme sensitive to the number of normal and tangential derivatives inside $\pa$. If we decompose $\pa$ in a general form
\be\label{E:DECOMP}
\pa = \pa_r^k\slashed\nabla^\beta, \ \ |\beta|+k = N,
\ee
then whenever we commute the PDE with an operator of the form~\eqref{E:DECOMP} we evaluate the inner product with $\pa\theta_\tau$ in the weighted $L^2_{w^{\alpha+k}}$ energy space, where
$\|f\|_{L^2_{w^{\alpha+k}}}^2=\int_{\Omega}f^2 w^{\alpha+k}\,dy$. Here $\slashed\nabla$ represents the tangential gradient. The role of the weight $w^{\alpha+k}$ is to offset the degeneracy caused by the vanishing of the enthalpy at the boundary. Using a high-order energy framework in such weighted spaces and Hardy-Sobolev embeddings one can adapt the well-posedness framework from~\cite{JaMa2015} to our setting. The damping effect mentioned above then allows us to extend our local solution to a global one and thus prove Theorem~\ref{T:SIDERISGLOBAL}. We refer the reader to~\cite{HaJa2016-2} for more details. The proof of Theorem~\ref{T:GWGLOBAL} follows similar ideas. It is different in that Goldreich-Weber stars can in fact expand at two different rates: the linear-in-$t$ rate (just like Sideris flows) and the non-generic self-similar rate $\sim_{t\to\infty} t^{\frac23}$. The latter solutions are only co-dimension $1$ stable and to prove this one necessitates a refined spectral analysis around the solutions of the generalised Lane-Emden equation~\eqref{E:GLE2}. We refer the reader to~\cite{HaJa2016-1} for the details  in the radial case. Very recently, nonradial stability 
in the class of irrotational perturbations has been shown in~\cite{HaJaLa2022}. 

\begin{remark}
Stability of the collapsing Goldreich-Weber stars is an open problem.
\end{remark}

\begin{remark}
The Sideris expanding flows in the regime $\ga>\frac53$ were shown to be stable by Shkoller and Sideris~\cite{ShSi}.
\end{remark}

\begin{remark}
Ideas related to the above expansion mechanism have been further developed in~\cite{P2022,PHJ2021,RHJ2021,Ri2021-1,Ri2021-2}.
\end{remark}

\subsection{Near-dust collapse}\label{SS:DUST}


\subsubsection{Self-gravitating dust}\label{SS:DUST}


A historically integral part of the mathematical study of stellar collapse is the collapse
of self-gravitating dust clouds. By definition, dust corresponds to a pressureless fluid, 
which amounts to the assumption $p=0$ in~\eqref{E:EPCONT}--\eqref{E:EPPOISSON}
or $p=0$ in~\eqref{E:EMTENSOR} in the relativistic case. The absence of pressure simplifies the problem 
considerably, to the point where one can explicitly solve the resulting PDE under the further assumption 
of radial symmetry. 

In the well-known works~\cite{To1934,Lemaitre1933} Tolman and Lema\^itre solved the 
dust-Einstein equations in radial symmetry, thus providing a large family of explicit solutions, parametrised
by the choice of the initial mass-energy density $\rho_0$. In their works it is very natural to use the 
comoving coordinates and we next briefly explain how to derive the Newtonian analogue 
of Tolman-Lema\^itre solutions using Lagrangian coordinates. We begin with a physically natural assumption that the initial density $\rho_0(r)$ is monotonically decreasing, i.e.
$\rho_0'<0$ on the interior of the support of $\rho_0$, which we take to be the interval $[0,1]$. \footnote{The monotonicity condition is for example satisfied by the Lane-Emden stars.}



If we assume that the dust-cloud is spherically symmetric and use the Lagrangian coordinates, then by~\eqref{E:EPRADIALLAGR}, the dust-Poisson system takes the form
\begin{align}\label{E:DUSTLAGR}
\chi_{tt}   + \frac{G(r)}{\chi^2} = 0,
\end{align}
where we recall~\eqref{E:AVERAGEDENSITY}.
This is remarkable, as in the Lagrangian variables, the problem reduces to a family of ODE parametrised by the choice of $r\in[0,1]$.
Given the initial conditions
\begin{align}\label{E:DUSTINITIAL}
\chi(0,r)=\chi_0(r)>0, \ \ \chi_r(0,r) = \chi_1(r),
\end{align}
the solution can thus be determined uniquely.
We consider the inward moving initial velocities with $\chi_1<0$. 
From the energy conservation 
\begin{align}\label{E:ENERGYDUST}
E(t) := \frac12\chi_t^2 - \frac{G(r)}{\chi} = E(0)
\end{align}
we obtain the formula
\begin{align}\label{E:PDEDUST2}
\chi_t= - \sqrt{\chi_1^2 + 2G\left(\frac1{\chi}-\frac1{\chi_0}\right)}.
\end{align}
Integrating~\eqref{E:PDEDUST2} one sees that for every $r$ there exists a $0<t^\ast(r)<\infty$ such that $\chi(t^\ast(r),r)=0$. 
A simple calculation reveals that for any $r\in[0,1]$ 
we have the universal blow-up exponent $2/3$ 
\begin{align}\label{E:UNIVERSALSCALING}
\chi(t,r) \sim c(r) (t^\ast(r) - t)^{\frac23}, \ \ t\to t^\ast(r). 
\end{align}

\begin{figure}

\begin{center}
\begin{tikzpicture}

\begin{scope}[scale=0.6, transform shape]
\coordinate [label=below:$0$] (A) at (-3,7){};
      \coordinate [label=below:$1$] (B) at (1.5,7) {};
      \coordinate [label=below:$r$] (C) at (-1,7) {};
      \coordinate [label=left:$t$] (E) at (-3,9) {};

\coordinate [label=left:$\frac1{g(0)}$] (D) at (-3,11){};
\coordinate  (L) at (0.9,13.0){};
\coordinate []  (M) at (1.5,14.6){};
\tkzCircumCenter(D,L,M)\tkzGetPoint{O}
\draw [fill=gray!20, very thick] (D) -- (A) -- (B) -- (M);
\draw[fill=white] (D) .. controls +(2,-0.5) and +(0.1,-0.6) .. (M);
\coordinate [label=below: $\mathcal G$] (B) at (-1,9.2) {};
\node at (-1,12) {$\Gamma$};
\end{scope}
\end{tikzpicture}
 \caption{Dust collapse in Lagrangian coordinates}
    \label{F:DUSTSUPPORT}
\end{center}
\end{figure}
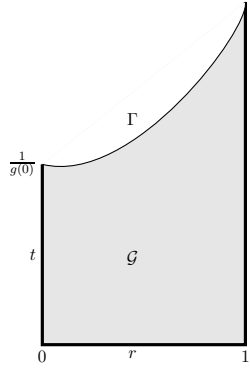

We may further define
the first blow-up time  
\[
t^\ast:=\min_{r\in[0,1]}t^\ast(r).
\]
Observe that the Eulerian description of the solution seizes to make sense at and after time $t\ge t^\ast$, as for example $\J[\chi]\big|_{t=t^\ast}=0$, see~\eqref{E:DETFORMULA}.
On the
other hand, for different values of $r$ the Lagrangian solution may make sense even after $t^\ast$. In particular, when $t^\ast(r)$ is a non-constant 
function, we can speak of the continued stellar collapse, wherein particles with a different Lagrangian label $r$ collapse at different times. 
This behaviour is generic and is referred to as inhomogeneous collapse~\cite{To1982}.
To further illustrate the nature of this type of singular behaviour we consider a particular subclass of solutions  given by 
\begin{align}\label{E:FUNDPROFILE0}
\chid(t,r) = (1 - g(r)t)^{\frac23}, \ \ g(r):= \left(\frac{9G(r)}{2}\right)^{\frac12}.
\end{align}
Notice that this solution satisfies~\eqref{E:DUSTLAGR} with the initial values $\chi_0(r)\equiv1$ and $\chi_1(r) = - \frac23 g(r)<0$. In this case
the average density $G(r)$ is given by
\be\label{E:CAPITALG}
G(r) = \frac1{r^3}\int_0^r 4\pi \rho_0(s) s^2\,ds,
\ee
which is a decreasing function since $\rho_0$ is decreasing.
It follows that $\chid$ becomes zero along the space-time curve 
\begin{align}\label{E:GAMMADEF}
\Gamma : = \{(t,r)\, | \, 1-g(r)t = 0\},
\end{align}
see Figure~\ref{F:DUSTSUPPORT}.
The solution is only well-defined 
in the region 
\[
\mathcal G:=\{(t,r)\, \big|\, 1-g(r)t>0\}.
\]

After a simple calculation using~\eqref{E:DETFORMULA} we have
\[
\J[\chid](t,r) = (1-g(r)t)^2\left(1-\frac23\frac{trg'(r)}{1-g(r)t}\right), \ \ (t,r)\in \mathcal G.
\]
In particular, $\chid$ and $\J[\chid]$ vanish along $\Gamma$ and therefore, since the Eulerian density  satisfies
\begin{align}\label{E:DUSTDENSITY}
\rho_{\text{dust}}(t,\chid(t,r)x) = \rho_0(r) \J[\chid](t,r)^{-1}, \ \ r=|x|, \ (t,r)\in \mathcal G,
\end{align}
the value of $\rho_{\text{dust}}(t,0)$ diverges to infinity at the first blow-up time $t^\ast:=\frac1{g(0)}$. In the region $\chid>0$, the Eulerian density $x\mapsto \rho_{\text{dust}}(t,x)$ is always well-defined away from the origin $x=0$.
Moreover for any $r\in[0,1]$
\[
\lim_{t\to\frac1{g(r)}} \rho_{\text{dust}}(t,\chid(t,r)x)=\infty.
\]
Since $r\mapsto g(r)$ is monotonically decreasing, 
particles that start out closer to the boundary of the star take longer to vanish into the singularity.

Finally, the support of the collapsing dust star shrinks to zero as $t\to\frac{1}{g(1)}$. This is clear, as the free boundary in the Eulerian description is at distance $\chid(t,1)=(1-g(1)t)^{\frac23}$ from the origin.
As $t\to\frac1{g(1)}$ the star concentrates with its mass completely absorbed at the origin:
\[
\lim_{t\to \frac1{g(1)}}\chid(t,1)=0.
\]

Therefore the time $t=\frac1{g(1)}$ has a natural interpretation as the end-point of star collapse for the dust example considered here.



\subsubsection{Near-dust collapse}\label{SS:NEARDUST}


It is natural to ask if there exist initial data that lead to finite time collapse that is qualitatively 
similar to the dust collapse described above, but in the presence of pressure. 
This is an a priori challenging problem, because pressure enters the equation at the highest order
of regularity and counteracts density concentration. An affirmative answer to this question is formulated
in the following theorem.

\begin{theorem}[\cite{GHJ2021a}]\label{T:NEARDUST}
Let $\gamma\in(1,\frac43)$. Then there exist initial data $\chi_0,\chi_1$ to~\eqref{E:EPRADIALLAGR}, a natural number
$n=n(\gamma)\in\mathbb N$ sufficiently large, and a compactly supported initial density profile $\rho_0:[0,1]\to\mathbb R_+$, satisfying
\begin{itemize}
\item $\rho_0>0$ on $[0,1)$ (density is positive);
\item $\rho_0(1)=0$ (vacuum at $r=1$);
\item There exists a constant $c>0$ such that
\be\label{E:FLATNESS}
\rho_0(r) = \rho_0(0) - c r^n + O_{r\to0^+}(r^{n+1}),  \ \ \text{ (the star is sufficiently flat at the centre)},
\ee
\end{itemize}  
so that the resulting Euler-Poisson flow leads to finite time collapse along the spacetime surface 
$\Gamma$ defined in~\eqref{E:GAMMADEF}, with $g(r)$ as in~\eqref{E:FUNDPROFILE0} and~\eqref{E:CAPITALG}.
\end{theorem} 

To build some intuition, for any $\ga\in(1,\frac43)$ consider a mass-invariant scaling of the EP-system given by the change of variables 
$
\rho \mapsto \lambda^{-3} \rho(s,y),  \ \
\mathbf{u}(t,x) \mapsto \lambda^{-\frac{1}{2}} \mathbf{u}(s,y),
\ \ 
\phi(t,x) \mapsto \lambda^{-1}\phi(s,y),
$
where $s = \frac{t}{\lambda^{\frac32}}, \ \ y = \frac{x}{\lambda}.$
The rescaled system is similar to the original EP-system~\eqref{E:EPCONT}--\eqref{E:EPPOISSON} and reads
\begin{align}
\pa_t\rho+\text{div}(\rho{\bf u}) & = 0 \label{E:EPCONTRESCALED}\\
\rho D_t {\bf u} +\kl \nabla p + \rho \nabla \phi & = 0, \ \ \kl : = \l^{4-3\ga}, \label{E:EPMOMENTUMRESCALED}\\
\Delta\phi & = 4\pi\rho.\label{E:EPPOISSONRESCALED}
\end{align}
If we pretend for a moment that $\l$ corresponds to some dynamically shrinking scale that 
leads to density blow-up, then we see from the rescaled momentum equation~\eqref{E:EPMOMENTUMRESCALED}
that the coefficient $\kl=\l^{4-3\ga}$ is small when $\l\ll1$ if and only if $1<\ga<\frac43$. This reflects the mass-supercritical nature of the 
problem and suggests a possibility of taming the pressure term in~\eqref{E:EPMOMENTUMRESCALED}. This intuition is also
deeply problematic, as the pressure term enters the equations at the top order of regularity and it is not clear
whether the above scaling can conspire with the pressure to mimic a dust-like collapse described in Section~\ref{SS:DUST}.

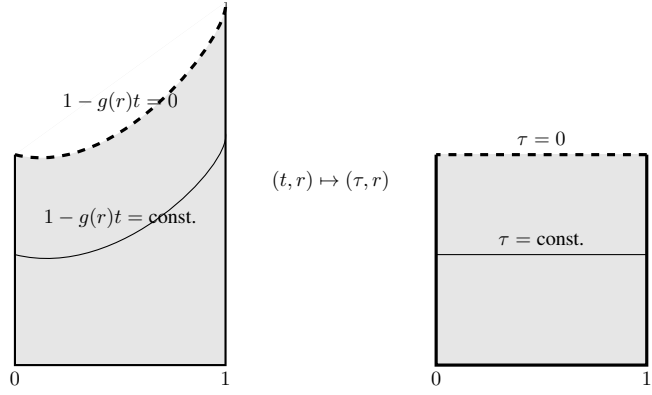
\begin{figure}
\begin{center}
\begin{tikzpicture}
\begin{scope}[scale=0.7, transform shape]
\coordinate [label=below:$0$] (A) at (-2,0);
      \coordinate [label=below:$1$] (B) at (2,0);

\tkzDefPoint(-2,1){P}\tkzDefPoint(-0.5,1.3){Q}\tkzDefPoint(2,2.4){S}

\tkzDefPoint(-2,2.1){A2}\tkzDefPoint(-0.6,2.6){M2}\tkzDefPoint(2,4.4){B2}

\tkzDefPoint(-2,4){D}\tkzDefPoint(1.2,5.5){L}\tkzDefPoint(2,6.9){M}

\draw [fill=gray!20, thick] (D) -- (A) -- (B) -- (M);

\tkzCircumCenter(A2,M2,B2)\tkzGetPoint{T2}

\draw[] (A2) .. controls +(2,-0.5) and +(0.1,-0.6) .. (B2);

\node at (0,2.8) {$1-g(r)t=\text{const.}$};

\tkzCircumCenter(D,L,M)\tkzGetPoint{O}
\draw[fill=white, dashed, very thick] (D) .. controls +(2,-0.5) and +(0.1,-0.6) .. (M);

\node at (0,5) {$1-g(r)t=0$};

\coordinate [label=below:$0$] (X) at (6,0){};
      \coordinate [label=below:$1$] (Y) at (10,0) {};
        \coordinate [label=right:] (Z) at (10,4) {};
      \coordinate [label=left:] (U) at (6,4) {};
\draw [fill=gray!20, very  thick] (U) -- (X) -- (Y) -- (Z);
\draw [dashed, very thick] (Z) -- (U);

\coordinate [label=left:] (X1) at (6,1){};
      \coordinate [label=right:] (X2) at (6,2.1) {};
\coordinate [label=left:] (Y1) at (10,1){};
      \coordinate [label=right:] (Y2) at (10,2.1) {};
\draw [] (X2) to (Y2);

\node at (8,2.4) {$\tau =\text{const.}$};
\node at (8,4.3) {$\tau=0$};

\node at (4,3.5) {$(t,r)\mapsto (\tau,r)$};
\end{scope}
\end{tikzpicture}
\caption{Foliation by the level sets of $\chid$}
\label{F:FOLIATION}
\end{center}
\end{figure}

Nevertheless, if we formally set $\kl=0$ in~\eqref{E:EPMOMENTUMRESCALED} we obtain the dust-Poisson system. This suggests 
the idea to try to build a collapsing solution to the EP-system as a ``perturbation" of $\chid$. To that end we first rectify the 
space-time support $G$ of the dust solution (see Figure~\ref{F:DUSTSUPPORT}) by foliating it with the level sets of the function $\chid$.
In other words, we introduce a new timelike coordinate $\tau$ and a new unknown $\phi$ 
\be\label{E:CHANGEOFVARIABLES}
\tau = 1-g(r)t, \ \ \phi(\tau,r) = \chi(t,r),
\ee
see Figure~\ref{F:FOLIATION}.
The dust solution $\chid$~\eqref{E:FUNDPROFILE0} now corresponds to 
\begin{align}
\phi_0(\tau,r) \equiv \tau^{\frac23}.
\end{align}

A simple change of variables implies that $\phi$ solves the PDE
\begin{align}\label{E:PHIEQN}
\pa_{\tau\tau}\phi + \frac{2}{9\phi^2} + \kl P[\phi] = 0,
\end{align}
where $P[\cdot]$ is the pull-back of the second order pressure operator, 
which in the new coordinates reads
\begin{align}
P[\phi] & = \frac{\phi^2}{g(r)w(r)^\alpha r^2} \Lambda \left(w(r)^{1+\alpha} J[\phi]^{-\ga}\right), \\
J[\phi] & = \phi^2(\phi+\Lambda\phi).
\end{align}
Here $\Lambda$ is the pull-back of the scaling operator $r\pa_r$ and reads
\[
\Lambda =  M_g  \pa_\tau + r\pa_r, \ \
M_g ( \tau, r) := (\tau-1)r\pa_r (\log g).
\]
The fundamental feature of the operator $\Lambda$ is that in the vicinity of the origin it takes 
the approximate form
\be\label{E:LAMBDAAPPROX}
\Lambda\approx r^n\pa_\tau + r\pa_r,
\ee
where we recall~\eqref{E:FLATNESS}. We see here that the first nontrivial term in the Taylor expansion of $\rho_0$
about the origin introduces a new effective scale in the problem in the vicinity of $r=0$.

The main idea of the proof is to construct a sufficiently accurate approximate solution $\phi_{\text{app}}$ to the problem 
so that 
\begin{align}\label{E:DECOMPOSITION}
\phi = \phi_{\text{app}} + \tau^m \frac Hr,
\end{align}
where
\be\label{E:APPROX}
\phi_{\text{app}} = \phi_0 + \kl\phi_1 + \dots +\kl^M \phi_M,
\ee
where $1\ll m,M\in\mathbb N$ are to be determined, and the remainder term $\frac{H}{r}$ is to be of order 1. We emphasise here that the 
all the difficulty in the problem is concentrated at the (new) spacetime origin $(\tau,r)=(0,0)$. We therefore want to focus on the regime where
both $0<r\ll1, 0<\tau\ll1$. 

To construct the approximate solution, we plug in the expansion~\eqref{E:APPROX} in~\eqref{E:PHIEQN} and formally expand in the powers of $\kl$.
This leads to
\begin{align}
\pa_{\tau\tau}\left(\phi_0+\kl\phi_1+\dots\right) + \frac2{9\left(\phi_0+\kl\phi_1+\dots\right)^2} + \kl P[\phi_0+\kl\phi_1+\dots] = 0.
\end{align}
Matching the powers of $\kl$ we obtain  a hierarchy of ODE of the form
\begin{align}\label{E:HIERARCHY}
\pa_{\tau\tau} \phi_{j+1} - \frac{4}{9\tau^2} \phi_{j+1} = f_{j+1}\left[\phi_0,\dots,\phi_j\right], \ \ j\ge0,
\end{align}
where $f_{j+1}$ are complicated combinatorial expressions that depend only on $\phi_0,\phi_1,\dots\phi_j$, but do not depend on $\phi_{j+1}$ - and therefore
act as source terms in~\eqref{E:HIERARCHY}. For such an expansion to be meaningful, we want $\phi_{j+1}$ to be
effectively smaller than the previous iterate $\phi_j$ in the near-collapse regime $\tau\ll1$. To illustrate this point, consider
the equation satisfied by $\phi_1$:
\begin{align}\label{E:PHIONEEQN}
\pa_{\tau\tau}\phi_1 - \frac{4}{9\tau^2} \phi_{1} = -P[\phi_0].
\end{align}
Recalling~\eqref{E:LAMBDAAPPROX} it is not hard to see that to the leading order
in the vicinity of $(0,0)$ the quantity $P[\phi_0]$ takes the form 
\begin{align}
P[\phi_0] & \approx \tau^{\frac13-2\ga} r^{n-2} \frac{1}{\left(1+\frac{r^n}{\tau}\right)^\ga} 
 = \tau^{\frac43-2\ga-\frac2n}\frac{\left(\frac{r^n}{\tau}\right)^{1-\frac2n}}{\left(1+\frac{r^n}{\tau}\right)^\ga}.
\end{align}
Going back to~\eqref{E:PHIONEEQN}, since it is a second order ODE it stands to reason that we can pick a solution $\phi_1$ so that it ``gains" two powers of $\tau$ with respect to 
the source term. In other words we solve~\eqref{E:PHIONEEQN}, so that the resulting solution near the singular point $(0,0)$ takes on the leading order form
\begin{align}
\phi_1 \approx \underbrace{\phi_0}_{=\tau^{\frac23}}  \tau^{\frac83-2\ga-\frac2n}\frac{\left(\frac{r^n}{\tau}\right)^{1-\frac2n}}{\left(1+\frac{r^n}{\tau}\right)^\ga}.
\end{align}
We let 
\begin{align}
\delta:=\frac83-2\ga-\frac2n = 2 \left(\frac43-\ga-\frac1n\right).
\end{align}
It is now clear that for any $1<\ga<\frac43$ there exists a choice of $n\in\mathbb N$ such that $\delta>0$. Moreover, for any $r\in[0,1]$ and $\tau>0$ the fraction
$\frac{\left(\frac{r^n}{\tau}\right)^{1-\frac2n}}{\left(1+\frac{r^n}{\tau}\right)^\ga}$ is always bounded from above uniformly, and we conclude that 
\be\label{E:GAIN}
|\phi_1|\lesssim \tau^\delta |\phi_0| \ll \phi_0,  \ \ 0<\tau\ll1.
\ee
This $\delta$-gain of smallness for the iterate $\phi_1$ is the fundamental feature of our strategy. It relies on both supercriticality of the problem, and the near-flatness assumption~\eqref{E:FLATNESS} with $n>\frac1{\frac43-\ga}$ chosen sufficiently large. The above calculation also illustrates the emergence and importance of the new characteristic 
spacetime scale $\frac{r^n}{\tau}$ in the region near the singularity $(r,\tau)=(0,0)$. The gain~\eqref{E:GAIN} can be shown to be propagated with every next iterate and spacetime derivatives. In~\cite{GHJ2021a} we show inductively that 
\begin{align}\label{E:GAIN2}
\lv\pa_\tau^m (r\pa_r)^\ell \phi_j\rv \lesssim_{m,\ell,j} \tau^{\frac23+j\delta-m}, \ \ \ell,j,m\in\mathbb N.
\end{align}

With~\eqref{E:GAIN2} in hand, we can now probe the ansatz~\eqref{E:DECOMPOSITION} by plugging it back into the original PDE~\eqref{E:PHIEQN}.
This will result in a new PDE for the remainder term $H$ with a source term which depends explicitly on $\phi_{\text{app}}$, see~\eqref{E:APPROX}.
This will result in a quasilinear wave equation of the schematic form
\begin{align}\label{E:HPDE}
&g^{00}[\phi]\pa_{\tau\tau}H + g^{01}[\phi]\pa_{\tau r}H + \left(\frac{M(M-1)}{\tau^2}-\frac2{3\tau^2}\right)H - \kl c[\phi] \frac1{w^\alpha} \pa_r\left(w^{1+\alpha}\frac1{r^2}\pa_r\left(r^2H\right)\right) \notag\\
& = \mathcal S_M[\phi_{\text{app}}].
\end{align}
The coefficients $g^{00}[\phi]$, $g^{01}[\phi]$, and $c[\phi]$ all exhibit singular/degenerate behaviour in the vicinity of $(0,0)$ and therefore require care. To prove the existence of a solution $H$ to~\eqref{E:HPDE} we therefore develop a high-order energy scheme, which is suitably weighted with power-like weights in $\tau$ (to handle the degeneracy at $\tau=0$)
and powers of $w(r)$ (to handle the vacuum degeneracy at $r=1$). Most importantly, we can choose the ``depth" $M$ of the approximate solution $\phi_{\text{app}}$~\eqref{E:APPROX} sufficiently large to show that the source term $\mathcal S_M[\phi_{\text{app}}]$ is sufficiently small near $(\tau,r)=(0,0)$. This scheme then leads to 
a global existence result for all $\tau\in(0,1]$. By the change of variables~\eqref{E:CHANGEOFVARIABLES} the data at $\tau=1$ correspond to initial data at $t=0$, which by construction lead to a dust-like collapsing solution, which to the leading order behaves like~\eqref{E:FUNDPROFILE0}.

\begin{remark}
It is not clear whether there exist open classes of data in a reasonable topology that exhibit a dust-like collapse.  
\end{remark}

\begin{remark}
We believe that the ideas used in the proof  of Theorem~\ref{T:NEARDUST} have a wider range of applicability to other quasilinear wave equations. Independently, related ideas were used in the semi-linear context by Cazenave, Martel, and Zhao~\cite{CaMaZh2019,CaMaZh2020}.
\end{remark}

\subsection{Self-similar collapse}\label{SS:NEWTONIANCOLLAPSE}


We now turn our attention to the construction of self-similar solutions 
that lead to finite time gravitational collapse from smooth data. We shall see that 
this point of view is extremely fruitful, as it will additionally lead to new insights in the
construction of collapsing relativistic flows in Section~\ref{S:NAKED}.

Based on the scaling invariance~\eqref{E:SCALING}, we define the self-similar variable
\begin{align}
y=\frac{r}{\sqrt{\k}(-t)^{2-\gamma}}
\end{align}
and look for solutions to~\eqref{E:EPCONT}--\eqref{E:EOS} of the form
\begin{align}
\rho(t,r) & =(-t)^{-2}\tilde\rho(y),\label{E:RHOSS}\\
u(t,r) & =\sqrt{\k}(-t)^{1-\gamma}\tilde u(y).\label{E:USS}
\end{align}
We introduce the 
relative (sometimes called the wind-) velocity
\begin{align}\label{E:OMEGADEF}
\om:=\frac{u+(2-\ga)y}{y},
\end{align}
and
plug in~\eqref{E:RHOSS}--\eqref{E:USS} into~\eqref{E:EPCONT}--\eqref{E:EPPOISSON}. This leads to the following system of
ordinary differential equations
\begin{align}\label{E:RHOEQN}
\rho'=&\,\frac{y\rho h(\rho,\om)}{G(y;\rho,\om)},\\
\om'=&\,\frac{4-3\ga-3\om}{y}-\frac{y\om h(\rho,\om)}{G(y;\rho,\om)}\label{E:OMEGAEQN},
\end{align}
where
\begin{align}
h(\rho,\om) &=2\om^2+(\ga-1)\om-\frac{4\pi\rho\om}{4-3\ga}+(\ga-1)(2-\ga), \\
G(y;\rho,\om)&=\ga\rho^{\ga-1}-y^2\om^2.
\end{align}

In the isothermal case $\ga=1$, system~\eqref{E:RHOEQN}--\eqref{E:OMEGAEQN} takes the form
\begin{align}\label{E:RHOEQNISO}
\rho'=&\,\frac{2y\rho\om(\om-2\pi \rho)}{1-y^2\om^2},\\
\om'=&\,\frac{1-3\om}{y}-\frac{y\om (\om-2\pi\rho)}{1-y^2\om^2}\label{E:OMEGAEQNISO},
\end{align}

We shall be interested in smooth $C^\infty$ solutions of~\eqref{E:RHOEQN}--\eqref{E:OMEGAEQN} in the mass-supercritical range $1\le\ga<\frac43$. 
This amounts to solving the non-autonomous system of ODE~\eqref{E:RHOEQN}--\eqref{E:OMEGAEQN} on the domain $[0,\infty)$. Any such solution
can be interpreted as an imploding star solution whose density blows up at the rate $\frac1{t^2}$ on approach to singularity along the cones of constant $y\equiv c$, i.e.
\[
\rho(t,x) = \frac1{t^2} \rho (c), \ \ \frac{|x|}{\sqrt \k (-t)^{2-\ga}} = c, \ \ t<0.
\]
In this formulation, the solution simply seizes to make sense as a classical solution at time $t=0^-$. Nevertheless, we shall see later that it makes sense to think of
$t=0^-$ as the first singularity and consider a singular curve in the $r-t$ plane along which the density blows up. A rigorous description of this
phenomenon necessitates the use of Lagrangian coordinates, see Section~\ref{S:NAKED}.

\begin{remark}
There is a parallel here to the dust and near-dust collapse from Sections~\ref{SS:DUST} and~\ref{SS:NEARDUST} in that they also allow for an extension 
that allows us to attribute different collapse times to different particles inside the gas. A decisive difference however, is that in the case of self-similar collapse, the rates of blow-up
at the scaling origin $\mathcal O=(0,0)$ (first time collapse point) and along the boundary of the maximal extension in the Lagrangian coordinates are different.
\end{remark}

The basic challenge in the analysis of~\eqref{E:RHOEQN}--\eqref{E:OMEGAEQN} is the possible presence of singularities in the denominator $G(y;\rho,\om)$ on the right-hand side. Any such point $y_\ast\in[0,\infty)$ by definition satisfies
\be\label{E:SONICDEF}
\ga \rho(y_\ast)^{\ga-1} - y_\ast^2\om(y_\ast)^2 = 0,
\ee
and is referred to as a sonic point. Geometrically it corresponds to the boundary of the backward acoustical cone emanating from the singularity at $(0,0)$.
Assuming smoothness at $y=0$ and that the solution is smooth and the density vanishes on approach to $y=\infty$, a simple Taylor expansion analysis
shows that we must have
\begin{align}
& \rho\sim_{y\to\infty} y^{-\frac{2}{2-\ga}}, \ \om\sim_{y\to\infty} 1, \label{E:FARFIELDLIKE}\\
& \rho(0) >0, \ \ \om(0) = \frac{4-3\ga}{3}.\label{E:FRIEDMANLIKE}
\end{align}


\begin{remark}\label{R:SPECIALSOLUTIONS}
There are in fact two explicit solutions to~\eqref{E:RHOEQN}--\eqref{E:OMEGAEQN}, each consistent with exactly one of the two boundary 
conditions~\eqref{E:FARFIELDLIKE}--\eqref{E:FRIEDMANLIKE}. They are the far-field solution
\begin{align}
\label{E:FARFIELDDEF}
\om_f=2-\ga,\quad \rho_f=ky^{-\frac{2}{2-\ga}},\text{ where }k=\Big(\frac{\ga(4-3\ga)}{2\pi(2-\ga)^2}\Big)^{\frac{1}{2-\ga}},
\end{align}
and the Friedman solution
\begin{align}\label{E:FRIEDMANPLP}
\om_F=\frac{4-3\ga}{3},\quad \rho_F=\frac{1}{6\pi}.
\end{align}
Note that $\rho_f$ blows up as $y\to0^+$ and $\rho_F$ does not decay to $0$ as $y\to\infty$. We may therefore, informally speaking, think of the smooth solution 
to~\eqref{E:RHOEQN}--\eqref{E:OMEGAEQN} satisfying~\eqref{E:FARFIELDLIKE}--\eqref{E:FRIEDMANLIKE} as a heteroclinic connection between the far-field and the Friedman solution.
\end{remark}

It is straightforward to see that any continuous solution $(\rho,\om)$ satisfying~\eqref{E:FARFIELDLIKE}--\eqref{E:FRIEDMANLIKE} leads to the boundary behaviour
$G(0;\rho,\om)>0$ and $\lim_{y\to\infty}G(y;\rho,\om)=-\infty$. Therefore by the intermediate value theorem any $C^1$-solution $(\rho,\om)$ has at least one sonic point $\bar y_\ast\in(0,\infty)$. So by necessity, we are interested in smooth solutions that feature a sonic point. We shall see that this is challenging, as the requirement of smoothness places severe restrictions on the behaviour of the solution at the sonic point.

\begin{remark}
In a recent pioneering work~\cite{MRRS1} Merle, Rapha\"el, Rodnianski, and Szeftel systematically constructed smooth self-similar solutions to the compressible Euler equations 
with the $\gamma$-law~\eqref{E:EOS} (for a.e. $\gamma>1$). Their work gives the first rigorous construction of physically relevant $C^\infty$ self-similar solutions featuring a sonic point in the setting of compressible fluids. 
The smoothness plays a crucial role in the dynamic stability analysis,  see~\cite{MRRS2}, wherein the finite co-dimension stability of the exact self-similar flows is shown. We emphasise that the underlying self-similar reduction gives an autonomous $2\times2$ system of ODE in contrast to  the non-autonomous nature of our system~\eqref{E:RHOEQN}--\eqref{E:OMEGAEQN}.  A numerical investigation of the space of imploding Euler flows was carried out by Biasi~\cite{Biasi}; for further existence results on smooth self-similar solutions of the Euler equation and their finite codimension stability, see a recent preprint~\cite{BuCLGS2022}.
\end{remark}

When $\gamma=1$ the $C^\infty$ solutions to the system~\eqref{E:RHOEQNISO}--\eqref{E:OMEGAEQNISO} were first analysed numerically in the seminal works of Larson~\cite{Larson1969} and Penston~\cite{Penston1969} in 1969. These works recognised the importance such a solution may play for the description of the late stage of the gravitational collapse of a Newtonian star. It is commonly referred to as the {\em Larson-Penston (LP)}-solution. Further milestone in the subject were the works of Hunter~\cite{Hunter77} and Shu~\cite{Shu77}, which among other things, suggested numerically the existence of a discrete family of $C^\infty$-self-similar solutions. In this picture, the LP-solution can be interpreted as the ``ground state" self-similar flow, which is expected to be the most stable member of the family. The next ``excited" state, the so-called Hunter A solution, plays an important role in the Newtonian criticality theory~\cite{MaedaHarada2001,GuGa2007}. Numerically, the LP-solutions are distinguished as the only member of the discrete family where the radial velocity $u$ is strictly positive. This implies that the flow is strictly inflowing in contrast to the excited solutions found by Hunter that exhibit regions of both positive and negative fluid velocity~\cite{Hunter77,MaedaHarada2001,BrWi1998}

Our main result in this setting is the existence of an LP-solution
\begin{theorem}[\cite{GHJ2021b}]\label{T:LP}
Let $\gamma=1$. Then there exists a real-analytic LP-solution $(\rho,\om)$ to~\eqref{E:RHOEQNISO}--\eqref{E:OMEGAEQNISO} satisfying the boundary conditions~\eqref{E:FARFIELDLIKE}--\eqref{E:FRIEDMANLIKE} with $\ga=1$.
\end{theorem}

Naturally, the assumption $\ga=1$ may seem restrictive.\footnote{Understanding the isothermal case $\ga=1$ is morally important for the fully relativistic setting, where the linear relation $p=\k\rho$ is indeed the correct equation of state at high densities~\cite{Ha1998,HaMa2001,HaLiRe2021}.} In fact, the physical relevance of the LP and Hunter solutions was criticised by Yahil~\cite{Yahil83} on the basis of having both infinite mass and energy. It is of course a separate question, whether one can truncate or suitably modify the LP solution in the asymptotic region so that we obtain a collapsing star of finite mass and energy (this is indeed possible, the proof is implicit in the results of~\cite{GHJ2021c}). But, if one insists on exact self-similarity then there can be no finite mass collapsing flows when $\gamma\in[1,\frac43)$. However, when $\gamma\in[\frac65,\frac43)$ self-similar solutions to~\eqref{E:RHOEQN}--\eqref{E:OMEGAEQN} do have finite energy, which can be read off from the way scaling affects the total energy~\eqref{E:ENERGYDEF}. This is the energy-critical ($\ga=\frac65$) and the energy-subcritical ($\ga\in(\frac65,\frac43)$) range in which Yahil~\cite{Yahil83} numerically established the existence of ground state analogues of the LP-solutions. We refer to them as Yahil solutions. Our next main result is the proof of existence of such solutions.

\begin{theorem}[\cite{GHJS2021}]\label{T:YAHIL}
For any $\gamma\in(1,\frac43)$ there exists a real-analytic Yahil-type solution $(\rho,\om)$ to~\eqref{E:RHOEQN}--\eqref{E:OMEGAEQN} satisfying the doundary conditions~\eqref{E:FARFIELDLIKE}--\eqref{E:FRIEDMANLIKE}. 
\end{theorem}

We devote the rest of this section to the main ideas behind the proofs of Theorems~\ref{T:LP}-\ref{T:YAHIL}. We will for expository brevity focus on the proof of Theorem~\ref{T:LP} and  indicate some conceptual differences to the proof of Theorem~\ref{T:YAHIL}. Theorem~\ref{T:YAHIL} is considerably more complex and its proof required some fundamental new ideas with respect to the proof of Theorem~\ref{T:LP}. 

By a slight abuse of notation we can rewrite~\eqref{E:RHOEQNISO}--\eqref{E:OMEGAEQNISO} in the form
\begin{align}\label{E:RHOEQNISO1}
\rho'=&\,\frac{2y\rho\om(\om-\rho)}{1-y^2\om^2},\\
\om'=&\,\frac{1-3\om}{y}-\frac{2y\om^2 (\om-\rho)}{1-y^2\om^2}\label{E:OMEGAEQNISO1},
\end{align}
where we changed the unknowns $\rho\mapsto 2\pi \rho$. The proof of Theorem~\ref{T:LP} contains three main steps.

{\em Step 1. Sonic point and the sonic window.}
As explained above, any smooth solution satisfying the asymptotic behaviour~\eqref{E:FARFIELDLIKE}--\eqref{E:FRIEDMANLIKE} must feature a sonic point, i.e. a point $y_\ast$ where $1-y_\ast^2\om(y_\ast)^2=0$. Such a singular point is a priori unknown and clearly depends on the solution itself. Our first task is to make an educated guess as to where one might expect such a point to be. It is not hard to see that the sonic point associated with the Friedman solution~\eqref{E:FRIEDMANPLP} in the case $\ga=1$ is precisely $y_\ast=3$. We choose $3$ to be the right end-point of our sonic window, while the left end-point is chosen to be $y_\ast=2$ for reasons that will become apparent below. Numerical investigations~\cite{BrWi1998,MaedaHarada2001} suggest that the sonic point value associated with the actual LP-solution is roughly $y_\ast\approx 2.41$.

Assuming the existence of a putative sonic point $y_\ast\in[2,3]$, our first goal is to show the existence of a smooth solution of the dynamical 
system~\eqref{E:RHOEQNISO1}--\eqref{E:OMEGAEQNISO1} in a small neighbourhood $(y_\ast-\delta,y_\ast+\delta)$ for some $\delta>0$. Assuming that there indeed exists a solution, it must have a Taylor expansion of the general form
\begin{align}\label{E:POWERSERIES}
\rho(y) = \sum_{N=0}^\infty \rho_N(y-y_\ast)^N, \ \ \omega(y) =  \sum_{N=0}^\infty \om_N(y-y_\ast)^N.
\end{align}
Our key idea is to prove convergence of the above series to obtain a local-around-$y_\ast$ real-analytic solution. Since the denominator $1-y_\ast^2\om(y_\ast)^2$ vanishes
we infer that i) $\om_0 = \frac1{y_\ast}$ and ii) as a consistency condition the numerator must vanish at $y_\ast$, i.e. $\om_0=\rho_0$. In conclusion
\[
\rho_0=\om_0=\frac1{y_\ast}.
\]
We then plug this information back into~\eqref{E:RHOEQNISO1}--\eqref{E:OMEGAEQNISO1} and expand to the next order in order to obtain the values of $\om_1$ and $\rho_1$. In doing so, an interesting dichotomy arises: it turns out that $\om_1$ satisfies a quadratic equation and two choices/types of solutions  emerge: $(\rho_1,\om_1)$ is
\begin{itemize}
\item either of Type 1: 
\be\label{E:TYPEONE}
\rho_1 = -\frac1{y_\ast^2}, \ \ \om_1 = -\frac1{y_\ast}(1-\frac2{y_\ast}),
\ee
\item or of Type 2:
\be\label{E:TYPETWO}
\rho_1 = \frac1{y_\ast}\left(1-\frac3{y_\ast}\right), \ \ \om_1 = 0.
\ee
\end{itemize}
It was recognised early~\cite{Hunter77} that the numerical solutions discovered by Larson and Penston are in fact of Type 1.   
Upon fixing the choice~\eqref{E:TYPEONE} it is then shown that for all $N\ge 2$ the Taylor coefficients $(\rho_N,\om_N)$  
are uniquely determined from the recursive relationship
\[
\mathcal A_N\begin{pmatrix} \rho_{N} \\ \om_N \end{pmatrix} = \begin{pmatrix} \mathcal F_{N} \\ \mathcal G_N \end{pmatrix},
\]
where $\mathcal F_N,\mathcal G_N$ are source terms that depend only on $(\rho_k,\om_k)$, $k\le N-1$ and
\begin{align}
\mathcal A_N(\omega_0,\omega_1,\rho_1) = 
\begin{pmatrix} -2N+2-2N\frac{\omega_1}{\omega_0} &  -\frac{2\rho_1}{\omega_0}-2 \\
-2 & -2N-4+\frac2{\omega_0}-(2N+2)\frac{\omega_1}{\omega_0} \end{pmatrix}.
\end{align}

\begin{remark}\label{R:HUNTER}
A similar recursion can be derived if we make the choice~\eqref{E:TYPETWO} instead. The choice~\eqref{E:TYPETWO} subsequently led Hunter to numerically discover a whole discrete 
family of further smooth self-similar solutions, sometimes labelled as Hunter A, B, C, D, \dots solutions, mentioned above. They are particularly interesting as they provide a link to so-called criticality theory in gravitational collapse~\cite{GuGa2007,NeCh}. It is an important open problem to construct Hunter solutions rigorously. They should be thought of as excited states in a discrete family of self-similar blow-up solutions, the ground state being the LP-solution. It is likely that 
the methods developed in~\cite{GHJ2021b,GHJS2021} will play an important role in such a construction.
\end{remark}

\begin{remark}
A surprising outcome of the sonic-point analysis in the polytropic setting $\ga>1$ is that the two seemingly disjoint branches~\eqref{E:TYPEONE}--\eqref{E:TYPETWO}
lend themselves to a different interpretation. They can be naturally viewed as members of the same branch which exhibits a discontinuity in the limit as $\ga\to1^+$.
\end{remark}

In~\cite{GHJ2021b} it is shown that the formal power series~\eqref{E:POWERSERIES} indeed converges on a small open interval around $y=y_\ast$. The challenge here is that 
there are no known general ODE results that can be cited to show the convergence in presence of sonic point singularities. Instead we resort to a combinatorial-type argument where we set up a careful inductive strategy and show that there exists a constant $C>0$ such that 
\begin{align}
|\rho_N|, |\om_N| \le \frac{C^N}{N^2}, \ \ N\ge 2.
\end{align}
This in turn implies the power series convergence. In conclusion, for any $y_\ast\in[2,3]$ we obtain a solution
\begin{align}
(\rho(\cdot;y_\ast),\om(\cdot;y_\ast))
\end{align}
defined on an interval $(y_\ast-\delta,y_\ast+\delta)$, where $\delta>0$ can be chosen independently of $y_\ast\in[2,3]$. By standard continuity arguments 
we extend such a solution to its maximal interval of existence to the left of the form $(s(y_\ast),y_\ast]$, where $s(y_\ast)\in[0,y_\ast-\delta]$ for all $y_\ast\in[2,3]$.

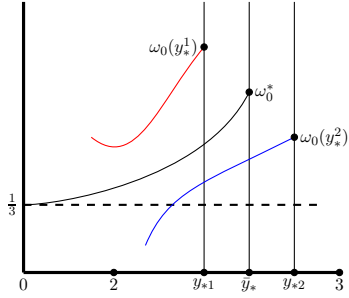
\begin{figure}
\begin{center}
\begin{tikzpicture}
\begin{scope}[scale=0.6, transform shape]
\coordinate[label=below:$0$] (A) at (0,0){};
\coordinate[label=below:$y_{\ast1}$] (B) at (4,0){};
\coordinate[label=below:$\bar y_\ast$] (C) at (5,0){};
\coordinate[label=below:$y_{\ast2}$] (D) at (6,0){};
\coordinate[label = left:$\frac13$] (E) at (0,1.5){};
\coordinate[label = below:$2$] (K) at (2,0){};
\coordinate[label = below:$3$] (L) at (7,0){};

\draw[fill=black] (K) circle (2pt);

\draw[fill=black] (L) circle (2pt);

\draw[fill=black] (B) circle (2pt);

\draw[fill=black] (C) circle (2pt);

\draw[fill=black] (D) circle (2pt);

\draw[very thick] (A)--(7,0);

\draw[very thick] (A)--(0,6);

\draw[dashed, thick] (E)--(6.5,1.5);

\draw (B) -- (4,6){};

\draw (C) -- (5,6){};

\draw (D) -- (6,6){};

\coordinate[label=right:$\om_0(y_\ast^2)$] (F) at (6,3){};

\draw[blue] (F) .. controls +(-2,-1) and +(0.5,1.2) .. (2.7,0.6);

\draw[fill=black] (F) circle (2pt);

\coordinate[label=right:$\om_0^\ast$] (G) at (5,4){};

\draw (G) .. controls +(-1,-2) and +(.7,0) .. (E);

\draw[fill=black] (G) circle (2pt);

\coordinate[label=left:$\om_0(y_\ast^1)$] (H) at (4,5){};

\draw[red] (H) .. controls +(-1,-1.2) and +(1,-0.8) .. (1.5,3);

\draw[fill=black] (H) circle (2pt);
\end{scope}
\end{tikzpicture}
 \caption{Schematic depiction of the shooting argument}
 \label{F:SHOOTING}
\end{center}
\end{figure}

{\em Step 2. The Friedman connection.}
Our goal is to identify a $\bar y_\ast\in[2,3]$ such that the local solution $(\rho(\cdot;y_\ast),\om(\cdot;y_\ast))$ exists all the way to $y=0$ and satisfies the boundary condition 
\be\label{E:FB}
\om(0;y_\ast)=\frac13.
\ee 
To this end we develop an effective shooting method, which turns out to be a very robust idea applicable also in the relativistic case, see Section~\ref{S:NAKED}. Using the precise knowledge of the Taylor coefficients in the expansion for $(\rho,\om)$ near the sonic point, one can show that for all $y_\ast<3$, but sufficiently close to $y_\ast=3$, the associated local LP-type solution ``dips" below $\om=\frac13$. Assuming that $y_{\frac13}(y_\ast)$ is the supremum of all such points - i.e. the ``first" point going to the left from the sonic point so that $\om(y_{\frac13}(y_\ast);y_\ast)=\frac13$, we see from~\eqref{E:OMEGAEQNISO1} that 
\be\label{E:OMEGA13}
\omega'(y_{\frac13}) = -\frac{2y_{\frac13}(\frac13-\rho(y_{\frac13}))}{9(1-\frac{y^2}9)}.
\ee
On the other hand, by subtracting~\eqref{E:RHOEQNISO1} from~\eqref{E:OMEGAEQNISO1} we also get
\begin{align}
(\om-\rho)' = \frac{1-3\om}{y} -\frac{2y\om (\om-\rho)}{1-y^2\om^2} \left(\rho+\om\right).
\end{align}
Based on this relation and the strict inequality $\rho-\om>0$ which holds ``initially" on a small interval $(y_\ast-\delta,y_\ast))$ to the left of the sonic point, we can show that
$\rho>\om$ on the maximal interval of existence (to the left) of the solution. It follows in particular that $\omega$ dips strictly below $\frac13$ locally to the left of $y_{\frac13}$. If $\omega$ were to ever cross $\frac13$ again on its maximal interval of existence, this would contradict the relation~\eqref{E:OMEGA13} and the strict inequality $\rho>\om$. Therefore
the region $\{y\,\big| \omega<\frac13\}$ is invariant under the flow, and as a consequence any solution that does go below $\frac13$ cannot be a regular solution of the problem as it cannot satisfy the boundary condition~\eqref{E:FB}. This then suggests a strategy to construct a regular connection. We define $\bar y_\ast$  via
\begin{align}
\bar y_\ast &: = \inf_{y_\ast\in Y} y_\ast,\label{E:BARYDEF} 
\end{align}
where
\begin{align}\label{E:YDEFA}
Y : = \left\{y_\ast\in[2,3]\,\big| \ \exists \, y\in (s(\tilde y_\ast), y_\ast) \ \text{such that } \omega(y;\tilde y_\ast)=\frac13 \ \text{ for all } \ \tilde y_\ast\in[y_\ast,3]\right\}.
\end{align}
In other words, $\bar y_\ast$ is the ``first" $y_\ast$ smaller than $3$ such that  the associated solution $\om(\cdot;\bar y_\ast)$ never takes on a value strictly smaller than $\frac13$
on its maximal interval of existence to the left  $(s(\bar y_\ast),\bar y_\ast]\subset [0,y_\ast]$. This ``shooting" method is illustrated in Figure~\ref{F:SHOOTING}. With a little bit of extra work it is possible to show that this minimality property embedded in the definition~\eqref{E:BARYDEF} guarantees that $s(\bar y_\ast)=0$, i.e. the solution exists on the whole interval $(0,\bar y_\ast]$. The final step is to show that the thus constructed solution $\om(\cdot;\bar y_\ast)$ satisfies the boundary condition~\eqref{E:FB} and does not, for example, blow up on approach to $y=0$. In the work~\cite{GHJ2021b} this is accomplished through a series of rather technical lemmas, but the procedure can be considerably simplified by following the strategy introduced in the polytropic context~\cite{GHJS2021}. One can in fact show that $\om(\cdot;\bar y_\ast)$ is also a monotonically increasing function (i.e. decreasing in value as we move to the left). The proof is delicate and uses both the minimality property~\eqref{E:BARYDEF}--\eqref{E:YDEFA} and nonlinear invariances in a nontrivial way, see~\cite{GHJS2021} for details.


{\em Step 3. The far-field connection.}
The final step in the proof of Theorem~\ref{T:LP} is to ensure that we can solve~\eqref{E:RHOEQNISO1}--\eqref{E:OMEGAEQNISO1} ``to the right of the sonic point", i.e. on the interval $[\bar y_\ast,\infty)$ so that the asymptotic behaviour~\eqref{E:FARFIELDLIKE} (with $\ga=1$) is honoured. This far-field connection turns out to be  easier to construct than the  Friedman connection from Step 2 above, and the main heuristic reason is that the flow is more stable to the right of the sonic point. We use basic invariances of the system to ensure that for any $y_\ast\in[2,3]$ the solution $\om(\cdot,y_\ast)$ stays trapped in the interval $(\frac13,1)$ for all $y>y_\ast$, as well as to ensure the inequality $\omega>\rho$ for all $y>y_\ast$ on the maximal interval of existence. With these two properties in place, it is not hard to see that the solution $(\rho(\cdot;y_\ast),\om(\cdot;y_\ast))$ then exists globally 
on the interval $[y_\ast,\infty)$ for all $y_\ast\in[2,3]$ and satisfies the boundary condition~\eqref{E:FARFIELDLIKE}. Choosing $y_\ast=\bar y_\ast$ constructed in Step 2 above, we finally obtain a globally defined real analytic solution on $[0,\infty)$ satisfying the correct boundary conditions. The resulting flow is the Larson-Penston solution claimed in Theorem~\ref{T:LP}.


\section{Self-gravitating relativistic gases}\label{S:NAKED}

\begin{figure}

\begin{center}
\begin{tikzpicture}
\begin{scope}[scale=0.8, transform shape]

\coordinate [] (A) at (0,0){};

\coordinate [] (A') at (0,-0.5);

\coordinate[label=left:$\mathcal O$] (O) at (0,2);

\coordinate (B) at (4,4);

\coordinate (B') at (4.5,3.5);

\coordinate (C) at (3,5);

\coordinate (E) at (1,1);

\coordinate (F) at (3,3);

\coordinate (F2) at (3.6,3.6);

\coordinate (G) at (2,4);

\coordinate (G2) at (2.6,4.6);

\draw[dashed,->] (F)--(G);

\draw[dashed,->] (F2)--(G2);

\draw [thick] (O) -- (A);

\draw (A) -- (B);

\draw[dashed] (B) -- (C);

\draw [thick] (O) -- (C);

\draw (O)--(E);

\node at (0.6,1.7) {$\mathcal N$};

\node at (3.8,4.7) {$\mathcal I^+$};

\node at (1.1, 3.8) {$\mathcal C\mathcal H^+$};

\node at (-0.5,1) {$r=0$};

\end{scope}
 \end{tikzpicture}
   \end{center}
    \caption{An incomplete future null-infinity.}  
\label{F:NS}
\end{figure}
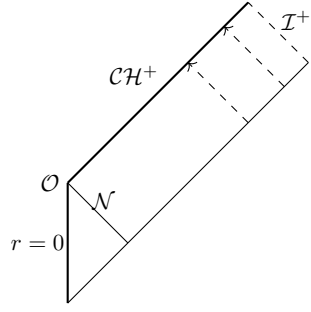

Our goal here is to explain the construction of a family of self-similar imploding relativistic solutions from~\cite{GHJ2021c}. These solutions are analogous to the Larson-Penston flows constructed above, but with a new twist. The analysis of the causal structure of such  {\em relativistic Larson-Penston} spacetimes implies that the spacetime contains a so-called  naked singularity. This notion on the other hand is intimately tied to the weak cosmic censorship hypothesis of Penrose~\cite{Penrose1969,Ch1999b,RoSR2019}. We work here with a rigorous mathematical definition of naked singularities given by Rodnianski and Shlapentokh-Rothman - Definition 1.1 in~\cite{RoSR2019}. In a nutshell, we say that a spacetime contains a naked singularity if 
\begin{itemize}
\item it is a maximal global hyperbolic development with a complete Cauchy hypersurface,
\item and it possesses an incomplete future null-infinity.
\end{itemize}
The latter by definition means
that exists a uniform upper bound on the affine length of a sequence of (suitably normalised) ingoing geodesics approaching null-infinity, starting from a regular outgoing cone - see Figure~\ref{F:NS}.

\subsection{Formulation of the problem and the main result}


Self-similar spacetimes $(M,g)$ are by definition characterised by the existence of a homothetic Killing vectorfield $\xi$ with the property
\[
\mathcal L_\xi g = 2g.
\]
We assume radial symmetry and wish to choose coordinates $(\tau,R)$ in a suitable conformal compactification so that $\xi=\tau\pa_\tau + R\pa_R$ indeed
corresponds to the usual scaling vectorfield. Following the pioneering work of Ori and Piran~\cite{OP1990}\footnote{Up to our knowledge, self-similar formulation of the EE-system was first given in~\cite{StShGu1965}, see also~\cite{CaTa1971}. For the derivation of the EE-system in comoving coordinates see for example~\cite{MiSh1964,EhKi1993}.} we work with the so-called comoving coordinates and the
metric takes the form
\be\label{E:METRIC}
g = -e^{2\mu(\tau, R)} d\tau^2 + e^{2\l(\tau,R)}dR^2 + r^2(\tau,R) \,\gamma.
\ee
Here $R$ is the comoving label, $\gamma$ the standard metric on $\mathbb S^2$, and the vectorfield $\pa_\tau$ is chosen to be parallel
to the $4$-velocity $u^\mu$, which in light of the normalisation~\eqref{E:NORMALISATION} implies $u=e^{-\mu}\pa_\tau$.
In this choice of coordinates, the Einstein-Euler system takes the form
\begin{align}
\pa_\tau\rho + (1+\k)\rho\left(\frac{\pa_R\mathcal V}{\pa_Rr}+2\frac {\mathcal V}r\right)e^\mu & =0,  \label{E:CONTEE}\\
\pa_\tau\l & = e^\mu \frac{\pa_R{\mathcal V}}{\pa_Rr},  \\
e^{-\mu}\pa_\tau \mathcal V + \frac \k{1+\k}\frac{\pa_Rr e^{-2\l}}{\rho} \pa_R\rho + 4\pi r \left(\frac13 G+\k\rho\right) & =0, \\
(\pa_Rr)^2 e^{-2\l} & = 1+ \mathcal V^2 - \frac{8\pi}{3} G r^2, \label{E:CONSTRAINTEE}
\end{align}
where
$
G(\tau,R): = \frac{m(\tau,R)}{\frac{4\pi}{3}r(\tau,R)^3}
$
is the average mass, 
$
\mathcal V : = e^{-\mu} \pa_\tau r
$
the radial velocity, and $\k$ the square of the speed of sound~\eqref{E:EOS}. In this formulation it is very easy to see that~\eqref{E:CONTEE}--\eqref{E:CONSTRAINTEE} is invariant under the scaling transformation
\begin{align*}
\rho\mapsto a^{-2}\rho(s,y), \ \ r\mapsto a r(s,y), \ \ \mathcal V \mapsto \mathcal V(s,y), \ \ \lambda\mapsto \lambda(s,y), \ \ \mu\mapsto \mu(s,y), 
\end{align*}
where the comoving ``time" $\tau$ and the label $R$ scale according to
\begin{align*}
 s = \frac{\tau}{a}, \ \ y = \frac{R}{a}, \ \ a>0. 
\end{align*}

\begin{figure}

\begin{center}
\begin{tikzpicture}
\begin{scope}[scale=0.6, transform shape]

\coordinate [label=left:$\mathcal O$] (A) at (-3,7){};

\coordinate [] (B) at (1,7){};

\node at (1,6.8) {$y=\infty$};

\node at (1,7.2) {$y=-\infty$};

\coordinate [] (C) at (-1.8,9.2){};

\coordinate[] (D) at (-2.5,4.5){};

\coordinate[] (E) at (0,5.8){};

\draw[dashed] (A)--(D);

\draw[dashed] (A)--(E);

\node at (0,5.7){$y=\text{const}$};

\node at (-2.3,4.6){$y=\text{const}$};


\coordinate [] (G) at (-3,9.3){};


\draw[very thick] (A)--(1,7);

\draw[very thick] (A)--(-3,4.5);

\draw[very thick] (A)--(-3,9);

\node at (-0.8,8.3) {$y<0$};

\node at (-0.8,5.2) {$y>0$};

\node at (-3.5,6) {$y=0$};

\draw[fill=white] (-3,7) circle (3pt);
\end{scope}
\end{tikzpicture}
\caption{Self-similar coordinates}
    \label{F:SSCOORD}
\end{center}
\end{figure}
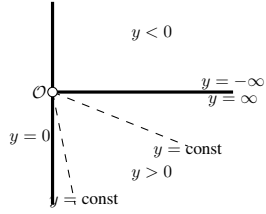

By analogy to the Newtonian case this motivates us to look for self-similar solutions, with the fundamental new variable
\[
y= \frac{R}{-\sqrt \k\tau}.
\]
We look for a solution to~\eqref{E:CONTEE}--\eqref{E:CONSTRAINTEE} of the form
\begin{align}
\rho(\tau,R) & = \frac{1}{2\pi \tau^2} \Sigma(y), \ \ r(\tau, R)  = -\sqrt \k \tau \tilde r(y),\ \
\mathcal V(\tau,R)  = \sqrt \k  V(y), \label{E:SS1}\\
\lambda(\tau,R) & = \lambda(y), \label{E:SS2}
\mu(\tau,R)  = \mu(y).
\end{align}

\begin{remark}
If we can indeed find a solution of the form~\eqref{E:SS1}--\eqref{E:SS2} such that  $\Sigma$ is bounded, then we conclude that $\rho(\tau,R)\asymp \frac1{\tau^2}$ as we approach the scaling origin $\mathcal O$ and we thus have implosion.
\end{remark}
Instead of working with the above variables directly, for the purpose of rigorous analysis, it is advantageous to work with 
the unknowns:
\begin{align*}
\d &: = \Sigma^{\frac{1-\k}{1+\k}}, \\
\w &: = (1+\k) \frac{e^\mu V + \tilde r }{\tilde r}  - \k . 
\end{align*}
It can be then checked that $(\d,\w)$ solve
\begin{align}
\d' &= -  \frac{ 2(1-\k) \d (\d-\w)}{(1+\k)y(e^{2\mu - 2\lambda} y^{-2} -1)}, \label{E:DEQN}\\
\w' &= \frac{(\w+\k)(1 -3\w)}{(1+\k)y} + \frac{2\w (\d-\w)}{y(e^{2\mu - 2\lambda} y^{-2} -1)}. \label{E:WEQN}
\end{align}

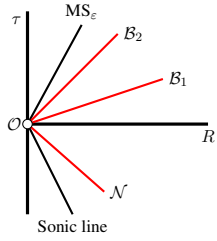
\begin{figure}
\begin{center}
\begin{tikzpicture}
\begin{scope}[scale=0.6, transform shape]

\coordinate [label=left:$\mathcal O$] (A) at (-3,7){};

\coordinate [label=below:$R$] (B) at (1,7){};

\coordinate [label=above:$\text{MS}_{\k}$] (C) at (-1.8,9.2){};

\coordinate [label=below:$\text{Sonic line}$] (D) at (-2,5){};

\coordinate [label=right:$\mathcal B_2$] (H) at (-1,9){};

\coordinate [label=right:$\mathcal B_1$] (E) at (0,8){};

\coordinate [label=right:$\mathcal N$] (F) at (-1.3,5.5){};

\coordinate [label=left:$\tau$] (G) at (-3,9.3){};


\draw[very thick] (A)--(1,7);

\draw[very thick] (A)--(-3,5);

\draw[very thick] (A)--(-3,9.5);

\draw[thick] (A)--(C);

\draw[thick] (A)--(D);

\draw[thick,red] (A)--(H);

\draw[thick,red] (A)--(E);

\draw[thick,red] (A)--(F);

\draw[fill=white] (-3,7) circle (3pt);
\end{scope}
\end{tikzpicture}
 \caption{Schematic depiction of the maximal self-similar extension
}
 \label{F:THMPIC1}
\end{center}
\end{figure}

\begin{remark}
The system~\eqref{E:DEQN}--\eqref{E:WEQN} is a non-autonomous $2\times2$ system and the sonic point corresponds to a point $y$ such that
\be
e^{2\mu - 2\lambda} y^{-2} -1=0.
\ee
The sonic point in fact corresponds to a hypersurface in the original spacetime $(\mathcal M,g)$, whose projection on the $(\tau,R)$-plane is schematically 
portrayed in Figure~\ref{F:THMPIC1}. Just like in the Newtonian case, it corresponds to the boundary of the backward acoustical cone emanating from the scaling origin $\mathcal O$.
\end{remark}

\begin{remark}
An approach inspired by dynamical systems ideas and phase space analysis to describe and classify smooth self-similar solutions of the EE-system was presented in~\cite{GoNiUg1998,CaCoGoNiUg2000}.
\end{remark}

The self-similar reduction~\eqref{E:DEQN}--\eqref{E:WEQN} originates from a manifestly Lagrangian description of the flow. On the other hand, the Newtonian version of the problem~\eqref{E:RHOEQNISO}--\eqref{E:OMEGAEQNISO} is derived from the Eulerian description. To link the two flows (for the moment just formally), we wish to formulate a suitable ``Eulerian" self-similar reduction of the problem. There is however no canonical choice of such coordinates in general relativity, but a choice that is well-suited to our needs are the so-called Schwarzschild coordinates, also used by Ori and Piran~\cite{OP1990}. At the level of self-similar reduction~\eqref{E:DEQN}--\eqref{E:WEQN} this amounts to introducing new variables
\[
x:=\tilde r(y),
\] 
\begin{align*}
\R(x):=\d(y), \ \ W(x):=\w (y).
\end{align*}
A calculation then reveals that 
\begin{align}
\R'(x) & = -  \frac{ 2x(1-\k) \R (W + \k) (\R-W) }{B},  \label{E:DEQNEULERIAN}\\
W'(x) & = \frac{(1 -3W )}{x} + \frac{2x(1+\k) W (W + \k) (\R-W)}{B},  \label{E:WEQNEULERIAN}
\end{align}
where 
\begin{align}\label{E:SONICDEN}
B=B[x;\R,W] : =   \R^{-\frac{2\k}{1-\k}} -\left[ ( W + \k)^2 - \k (W - 1)^2 + 4\k \R W \right] x^2,
\end{align}
Notice that in this description  sonic points correspond to the $x\in\mathbb R$ such that $B(x)=0$.
We make here an important observation - the formal
Newtonian limiting system is obtained by setting  $\k=0$, which corresponds precisely to~\eqref{E:RHOEQNISO}--\eqref{E:OMEGAEQNISO}.
This suggests that in the $0<\k\ll1$ regime we might hope to solve~\eqref{E:DEQNEULERIAN}--\eqref{E:WEQNEULERIAN} by emulating the 
strategy developed in Section~\ref{SS:NEWTONIANCOLLAPSE}. In that spirit first that by analogy to Remark~\ref{R:SPECIALSOLUTIONS} observe that there are two explicit solutions.
The far-field solution is given by 
\begin{align}\label{E:FARFIELDREL}
\R_f(x) = 
(1-\k)^{-\frac{2(1-\k)}{1+\k}} x^{-\frac{2(1-\k)}{1+\k}}, \ \
W_f(x)  = 1;
\end{align}
we note that the density $\R_f$  blows up at $x=0$ and decays to $0$ as $x\to\infty$. 
By contrast, the 
Friedmann
solution
\begin{align*}
\R_F(x)  = \frac13, \ \
W_F(x)  = \frac13,
\end{align*}
has bounded density at $0$, but it does not decay as $x\to\infty$.
Our goal is to find
a smooth solution connecting Friedmann to far-field, so that in particular
\[
\lim_{x\to\infty} W(x) =1, \ \ \lim_{x\to0^+}W(x) = \frac13.
\]
It is easily checked that this requirement implies the existence of at least one sonic point.

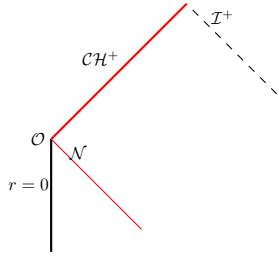
\begin{figure}

\begin{center}
\begin{tikzpicture}
\begin{scope}[scale=0.6, transform shape]

\coordinate [] (A) at (0,0){};

\coordinate [] (A') at (0,-0.5);

\coordinate[label=left:$\mathcal O$] (O) at (0,2);

\coordinate (B) at (4,4);

\coordinate (B') at (5,3);

\coordinate (C) at (3,5);

\coordinate (E) at (2,0);

\coordinate (F) at (0.5,.5);

\draw [thick] (O) -- (A');


\draw[dashed] (B') -- (C);

\draw [thick,red] (O) -- (C);

\draw[red] (O)--(E);


\node at (0.6,1.7) {$\mathcal N$};

\node at (3.8,4.7) {$\mathcal I^+$};

\node at (1.1, 3.8) {$\mathcal C\mathcal H^+$};

\node at (-0.5,1) {$r=0$};

\end{scope}
 \end{tikzpicture}
   \end{center}
\caption{Asymptotic flattening and naked singularities}
\label{F:THMPIC2}
\end{figure}

\begin{theorem}[~\cite{GHJ2021c}, Informal statement]\label{T:NS}
There exists an $0<\k_0<1$ such that for all $0<\k\le\k_0$ there exists a relativistic Larson-Penston solution (from now on RLP solution) with the following properties:
\begin{enumerate}
\item The maximal analytic extension of the self-similar spacetime has a singular boundary consisting of two components: the scaling origin  $\mathcal O=(0,0)$ and the Massive Singularity $\text{MS}_{\k}$ depicted in Figure~\ref{F:THMPIC1}. On approach to the scaling origin the fluid density / Ricci scalar blow up like $\tau^{-2}$;
\item There exists exactly one sonic line - the boundary of the backward acoustical cone emanating from $\mathcal O$, and exactly one null-curve $\mathcal N$ corresponding to the boundary of the backward null cone emanating from $\mathcal O$;
\item There exist at least two future oriented null-geodesics $\mathcal B_1$ and $\mathcal B_2$ emanating from the scaling origin $\mathcal O$ with the property that $r$ increases to $\infty$ as the affine parameter tends to $\infty$. Moreover, they are similarity curves of the spacetime (i.e. straight lines in the $R-\tau$ diagram) and $\mathcal B_1$ is the ``first" such outgoing null geodesic in the sense that any outgoing null-geodesic initiated in the causal past of $\mathcal B_1$ does not cross $\mathcal B_1$ nor it goes to/out of $\mathcal O$, see Figure~\ref{F:THMPIC1}.
\end{enumerate} 
Moreover, the solution can be flattened in the asymptotic region to give examples of naked singularity spacetimes. More precisely, for any $0<\k\le\k_0$ there exists a family of solutions to the Einstein-Euler system with the Penrose diagram like in Figure~\ref{F:THMPIC2}. The solution coincides with the exact self-similar RLP spacetime in the causal past of the curve $\mathcal N$, is asymptotically flat, and possesses an incomplete future null infinity. The singularity at $\mathcal O$ is visible to distant observers and the Cauchy Horizon $\mathcal C\mathcal H^+$ corresponds to the future null-geodesic emanating from $\mathcal O$ as in Figure~\ref{F:THMPIC2}, which effectively terminates the future null-infinity $\mathcal I^+$.
\end{theorem}


\begin{remark}
Our proof is nonperturbative - it is not the case that we use for example an implicit function theorem to perturb away from the Newtonian solution and construct the relativistic one. The nonlinear structure of the relativistic problem is rather delicate and does not lend itself to an easy perturbation argument. Instead, while our results do use the smallness of $\k$ in a crucial way, we do so to enforce the validity of various nonlinear invariances of the ODE flow, catered to the problem.
\end{remark}


\subsection{Strategy of the proof and main ideas}

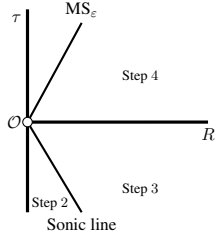
\begin{figure}

\begin{center}
\begin{tikzpicture}
\begin{scope}[scale=0.6, transform shape]

\coordinate [label=left:$\mathcal O$] (A) at (-3,7){};

\coordinate [label=below:$R$] (B) at (1,7){};

\coordinate [label=above:$\text{MS}_{\k}$] (C) at (-1.8,9.2){};

\coordinate [label=below:$\text{Sonic line}$] (D) at (-1.8,5){};



\coordinate [label=left:$\tau$] (G) at (-3,9.3){};


\draw[very thick] (A)--(1,7);

\draw[very thick] (A)--(-3,5);

\draw[very thick] (A)--(-3,9.5);

\draw[thick] (A)--(C);

\draw[thick] (A)--(D);

\node at (-2.5,5.2) {\footnotesize Step $2$};

\node at (-0.5,5.5) {\footnotesize Step $3$};

\node at (-0.5,8.0) {\footnotesize Step $4$};

\draw[fill=white] (-3,7) circle (3pt);
\end{scope}
\end{tikzpicture}
 \caption{Proof strategy} 
    \label{F:STRATEGY}
\end{center}
\end{figure}

The proof of Theorem~\ref{T:NS} proceeds in six steps.  The first four steps are devoted to the construction of the maximally extended RLP spacetime, in step 5 we study the causal structure of the self-similar spacetime, and finally in step 6 we explain how to flatten the self-similar solution asymptotically to obtain a spacetime with naked singularity in line with Figure~\ref{F:THMPIC2}.  Structurally, steps 1-3 echo the construction of the Larson-Penston spacetime explained in Section~\ref{SS:NEWTONIANCOLLAPSE} - step 2 corresponds to the relativistic version of the Friedman shooting argument, while step 3 establishes the existence of the far-field connection, see Figure~\ref{F:STRATEGY}. From the technical point of view, both steps are considerably more challenging than in the Newtonian case due to the highly nonlinear structure of the right-hand side of~\eqref{E:DEQNEULERIAN}--\eqref{E:WEQNEULERIAN}. Step 4 is new by comparison to the Larson-Penston analysis, as it provides an extension of the solution into the upper-half plane (i.e. positive times), which by definition is only possible in a suitable choice of comoving coordinates (i.e. not in the ``Eulerian" description). This leads to a nonstandard ODE problem and the question of how to extend the solution uniquely. In step 5 we demonstrate the existence of outgoing null-geodesics ${\mathcal B}_1$ and $\mathcal B_2$ and in Step 6 we show that $\mathcal B_1$ can be suitably deformed in the asymptotic region to correspond to the Cauchy horizon of an asymptotically flat spacetime with an incomplete future null-infinity in the sense of Christodoulou and Rodnianski--Shlapentokh-Rothman~\cite{Ch1999a,RoSR2019}.

{\em Step 1.} Sonic point $x_\ast$ is an a priori unknown singular point for the flow. By analogy to Step 1 of the proof of Theorem~\ref{T:LP}, we make 
an educated guess $x_\ast\in[2,3]$. This guess turns out to be slightly too optimistic and one instead works on the interval $[\xmin(\k),\xmax(\k)]$ where $\lim_{\k\to0}[\xmin(\k),\xmax(\k)]=[2,3]$, which is consistent with the Newtonian limit. The exact formulas for $\xmin(\k),\xmax(\k)$ are not so important - they play the same role as $2$ and $3$ do in the Newtonian problem in ensuring later on that our shooting argument works. We prove the local 
existence on an interval of the form $(x_\ast-\delta,x_\ast+\delta)$, where $\delta>0$ is independent $x_\ast$. We do this by considering the Taylor 
series expansion
\[
D = \sum D_N(x-x_\ast)^N, \ \ W = \sum W_N(x-x_\ast)^N
\]
at the sonic point $x_\ast$, anticipating that the requirement of smoothness will introduce a hierarchy of constraints for 
the sequence $(D_N,W_N)$, $N\ge0$. It is immediate from~\eqref{E:DEQNEULERIAN}--\eqref{E:WEQNEULERIAN} that 
we must have $W_0 = D_0$. By analogy to the Newtonian case, the next order coefficient $W_1$ is shown to satisfy
a quadratic equation with $W_0$-dependent coefficients. Much like in the Newtonian case, the two choices
are the relativistic analogues of choosing a Larson-Penston type or a Hunter-type expansion. In either case, 
the resulting sequence is uniquely specified once this choice is made, and we obtain a recursive relation of the form
\item
\[
A_N\begin{pmatrix} D_{N} \\ W_N \end{pmatrix} =  \begin{pmatrix} S^1_{N} \\ S^2_N \end{pmatrix}, \ \ N\ge 2,
\]
where
$S^j_N$ are source terms that depend only on $(D_k,W_k)$, $k\le N-1$ and $A_N$ is an invertible matrix whose entries depend explicitly on $x_\ast$ and $N\in\mathbb N$.

\begin{remark}
Note that the existence of Hunter-type solutions is not pursued in our work~\cite{GHJ2021c}. By analogy to the Newtonian case, see Remark~\ref{R:HUNTER}, the Hunter family of solutions provides a link to the criticality theory for relativistic gases~\cite{EvCo,CaCoGoNiUg2000, Ha1998,Ha2003,GuGa2007}, which has its roots in the well-known numerical work of Choptuik~\cite{Chop1993} on the collapse of massive scalar field. 
\end{remark}

{\em Steps 2 and 3: Friedman and far-field connections.}
Inspired by our Newtonian approach we wish to develop an effective shooting argument which will allow us to select
a value $\bar x_\ast$ so that the local solution around $\bar x_\ast$ of RLP-type extends all the way to the left to $x=0$.
We refer to it as the Friedman connection as any such solution, if smooth, must attain the value $W(0)=\frac13$. Similarly, in order to extend the solution to 
the right of the sonic point all the way to $x=\infty$ it must qualitatively behave like the far-field solution~\eqref{E:FARFIELDREL} - hence far-field connection. 
The proof to the left relies on 
a shooting procedure analogous to the one depicted in Figure~\ref{F:SHOOTING} which allows us to single out a specific value of $\bar x_\ast$ with the right properties. The analysis to the right is less sensitive to the choice of $x_\ast$ and indeed we show that for any value of $x_\ast\in[\xmin(\k),\xmax(\k)]$ there exists a global solution on $[x_\ast,\infty)$. 
However, in the relativistic case it is considerably harder to show that
the sign of the sonic denominator $B$ remains strictly positive (resp. negative) for all $x\in(0,\bar x_\ast)$ (resp. $x\in(x_\ast,\infty)$), by comparison to the Newtonian case.  A key new technical tool in our analysis is a Monotonicity Lemma designed to capture a novel nonlinear invariance in the problem, in the regime where $0<\k\ll1$.

We first recall the sonic denominator~\eqref{E:SONICDEN}
\[
B= \R^{-\frac{2\k}{1-\k}} -\left[ ( W + \k)^2 - \k (W - 1)^2 + 4\k \R W \right] x^2.
\]
To control its sign we introduce the following simple algebraic 
decomposition:
\[
B[x;\R,W] =  \left[\F -xW(x)\right] \left[H+ xW(x)\right],
\]
where
\begin{align*}
\F & : = -\frac{2\k}{1-\k}(1 +\R)x + \sqrt{ \frac{4\k^2}{(1-\k)^2} (1 +\R)^2 x^2+  \k x^2 
+ \frac{\R^{-\frac{2\k}{1-\k}}}{1-\k}} ,  \\
H& : =  \F +\frac{4\k}{1-\k}(1+\R)x. 
\end{align*}
It is now clear that the sign of 
$B$ corresponds to the sign of $\F-xW$. The key quantity to monitor turns out to be
\be
f(x) : = \F - x\R,
\ee
and a direct calculation~\cite{GHJ2021c} shows
that for any smooth solution of~\eqref{E:DEQNEULERIAN}--\eqref{E:WEQNEULERIAN} we have
\begin{lemma}[Monotonicity Lemma]\label{L:MON}
\[
f' + a[x;\R,W] f = b[x;\R,W],
\]
where $b$ has ``good sign properties in the regime we care about".
\end{lemma}
To show that the sign of $B$ remains positive to the left of the sonic point, we then show with the help of Lenma~\ref{L:MON} that the following ``sandwich" 
relations
\begin{align}
xW<&x\R<\F, \ \ \ x<x_\ast \\
xW>&x\R>\F, \ \ \ x>x_\ast
\end{align}
are dynamically propagated to the left (resp. right) of $x=\bar x_\ast$. It then remains to show that it is also satisfied in a small left (resp. right) neighbourhood of $x=\bar x_\ast$ which is accomplished using the precise information on the Taylor coefficients obtained in Step 1.

{\em Step 4. Maximal analytic extension.}
The metric $g$ becomes singular on approach to the horizontal axis $\tau=0$, i.e. as $x\to\infty$ or equivalently $y\to\infty$. We prove this by using the precise  
asymptotic behaviour as $x\to\infty$ and show that 
\begin{align}\label{E:MUAS}
e^{2\mu}\asymp_{y\to\infty}  y^{\frac{4\k}{1+\k}},
\end{align}
where we have conveniently reverted back to the Lagrangian variable $y$.
The goal is to show that the singular behaviour~\eqref{E:MUAS} is merely a coordinate singularity and it is possible to extend the metric smoothly (in fact analytically) across the surface $\tau=0$ into the upper half plane. We achieve this by introducing adapted comoving coordinates:
\begin{align}\label{E:BIGYDEF}
Y = y^{-\frac{1-\k}{1+\k}}, 
\end{align}
and a new set of unknowns
\begin{align}
\chi(Y): = \frac{\tilde r(y)}{y}, \ \ \RY(Y) : =  \d(y),  \ \ w(Y) : = \wl (y).
\end{align}
The choice~\eqref{E:BIGYDEF} of $Y$ is clearly motivated by~\eqref{E:MUAS}. A simple calculation reveals that the
system~\eqref{E:DEQN}--\eqref{E:WEQN} in the new variables takes the form
\begin{align}
\RY' &=    \frac{2\chi^2}{Y} \frac{ \RY (\vY+\k)^2(\RY-\vY)}{\mathcal C}, \label{E:NEW1}\\
\vY' &= -\frac{(\vY+\k)(1-3\vY)}{(1-\k)Y} -\frac{2(1+\k)\chi^2}{(1-\k)Y} \frac{\vY(\vY+\k)^2(\RY-\vY)}{\mathcal C}, \\
\chi' & = \frac{1-\vY}{(1-\k)Y} \chi,\label{E:NEW3}
\end{align}
where
\begin{align*}
\mathcal C : = \left(\RY Y^{-2}\right)^{-\eta}Y^2 - \chi^2 \left[(\vY + \k)^2-\k(\vY-1)^2 + 4\k \vY \RY \right]
\end{align*}
is the sonic denominator in the new variables,
The sharp asymptotic understanding 
of the solution as $y\to+\infty$ allows us to feed ``initial data" for the new system of ODE at $Y=0$ ($y=\infty$)
and show a local-in-$Y$ well-posedness result around $Y=0$. However, to extend the solution to the maximal interval of existence to the left of $Y$ we then prove
the following key lemma.

\begin{lemma}[The sandwich lemma] \label{L:SANDWICH}
On the maximal interval of existence to the left of $Y=0$,
we have the dynamic ``sandwich" bound:  
\[
1\lesssim \frac{\RY}{\vY} <1, \ \text{ for } Y<0.
\]
\end{lemma}

\begin{figure}
\begin{center}
\begin{tikzpicture}
\begin{scope}[scale=0.4, transform shape]

\coordinate [label=below:$0$] (A) at (-3,7){};
\coordinate [label=below:$\yms$] (B) at (-7,7){};
\coordinate [] (C) at (-7,11){};
\coordinate [label=below:$Y$] (D) at (1,7){};
\coordinate [label=right:$1$] (F) at (-3,9.2){};
\coordinate [label=above:$\vY$] (F) at (-5,10.3){};

\draw[very thick] (-3,7)--(-8,7);
\draw[very thick] (A)--(1,7);
\draw[very thick] (A)--(-3,11);
\draw[dashed, thick] (B) -- (-7,13);
\draw (-6.8,14) .. controls +(0.8,-4) and +(-1,0.5) .. (-3,9);
\draw (-3,9) .. controls +(1,-0.5) and +(-1,0) .. (1,7.8);

\draw[fill=black] (-3,7) circle (2pt);
\draw[fill=black] (-7,7) circle (2pt);
\draw[fill=black] (-3,9) circle (2pt);

\end{scope}
\end{tikzpicture}
\caption{Schematic depiction of $\vY(Y)$. Note that $\vY(Y)$ in the maximal extension blows up on approach to $\yms$}
\label{F:NEWODE}
\end{center}
\end{figure}
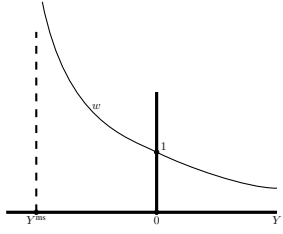

The proof of Lemma~\ref{L:SANDWICH} relies on careful exploitation of nonlinear invariances of~\eqref{E:NEW1}--\eqref{E:NEW3}. This in turn allows us to show that the maximal interval of existence is characterised by the simultaneous blow up of $d,w,$ and $\frac1\chi$, see Figure~\ref{F:NEWODE}, which is proved by showing that $w$ solves a Riccati type equation.

{\em Step 5. Causal structure and the existence of an outgoing null geodesic from $\mathcal O$.}

\begin{figure}

\begin{center}
\begin{tikzpicture}
\begin{scope}[scale=0.6, transform shape]

\coordinate [label=left:$\mathcal O$] (A) at (-3,7){};

\coordinate [] (B) at (1,7){};

\node at (1,6.8) {$y=\infty$};

\node at (1,7.2) {$y=-\infty$};

\coordinate [] (C) at (-1.8,9.2){};

\node at (-1.7,9.3) {$y=y_{\text{ms}}$};

\coordinate [] (G) at (-3,9.3){};


\draw[very thick] (A)--(1,7);

\draw[very thick] (A)--(-3,5);

\draw[very thick] (A)--(-3,9.5);

\draw[thick] (A)--(C);

\node at (-0.8,8.3) {$y<0$};

\node at (-3.5,6) {$y=0$};

\draw[fill=white] (-3,7) circle (3pt);
\end{scope}
\end{tikzpicture}
\caption{Set-up for the intermediate value theorem}
    \label{F:MVT}
\end{center}
\end{figure}
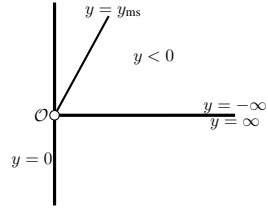

To complete the proof of Theorem~\ref{T:NS} we must show that there exists an outgoing radial null geodesic (RNG), denoted by $\mathcal B_1$ in Figure~\ref{F:THMPIC1},
emanating from $\mathcal O$ and reaching  infinity. To that end we recall the metric expression~\eqref{E:METRIC} - an RNG by definition satisfies $(\frac{dR}{d\tau})^2= e^{2\mu-2\l}$. Here we switched back to the comoving coordinates $(\mu,\l,y)$ which continue to make sense in the upper half-plane strictly away from the horizontal axis $\tau=0$. We first restrict our attention to a particular subclass of such geodesics that are simultaneously
 similarity lines in the problem so that there exists a $\sigma\in\mathbb R$ such that  $R(\tau) = \sigma \tau$ (we refer to them as simple RNG-s following the terminology of Ori and Piran~\cite{OP1990}). This trivially implies the relation
\[
y(\tau)=\frac{R(\tau)}{-\sqrt\k \tau}=-\frac{\sigma}{\sqrt\k},
\]
so we conclude that $\frac{\sigma}{\sqrt\k}$ is a root of
\begin{align*}
& y \mapsto F_\k (y) \\ 
& F_\k(y) : = \k y^2 e^{2\l-2\mu}-1
\end{align*}
A simple, but effective idea to prove the existence of a root of $F_\k$ is to use intermediate value theorem and goes back to Ori and Piran~\cite{OP1990}. The sharp asymptotic description of the self-similar flow from Steps 3 and 4 allows us to show that
\begin{itemize}
\item 
\be\label{E:YMSBEHAVIOUR}
\lim_{y\to y_\text{ms}} F_\k(y) = \lim_{y\to -\infty} F_\k(y) = \infty
\ee
\item 
and there exists a $y_0\in(-\infty,y_{\text{ms}})$ such that
\be\label{E:DIPPOINT}
F_\k(y_0)<0
\ee
if $\k$ is sufficiently small.
\end{itemize}
By the intermediate value theorem there must exist at least two 
slopes $y_1$ and $y_2$ corresponding to outgoing null-geodesics, which correspond precisely to curves $\mathcal B_1$ and $\mathcal B_2$ on Figure~\ref{F:THMPIC1}.
We emphasise here that the proof of~\eqref{E:YMSBEHAVIOUR} relies on the leading order asymptotics of the dynamical system~\eqref{E:NEW1}--\eqref{E:NEW3} on approach to the terminal point $Y=Y_{\text{ms}}$, or equivalently $y=y_{\text{ms}}$, while the second limit claimed in~\eqref{E:YMSBEHAVIOUR} mirrors the coordinate singularity at the horizontal axis
$\{\tau=0\}$ explained in Step 4. The key is therefore to establish~\eqref{E:DIPPOINT}, which schematically speaking, follows from very precise $\k$-independent bounds on the solution of~\eqref{E:NEW1}--\eqref{E:NEW3} in a small (but $\k$-independent) left neighbourhood of $Y=0$ (or equivalently $y=-\infty$), see Figure~\ref{F:MVT}.
The full description of (not necessarily simple) RNG-s as well as nonradial null-geodesics can be found in Section 8 and Appendix B of~\cite{GHJ2021c}. We note that after the publication~\cite{OP1990}, nonradial null-geodesics were analysed by Joshi and Dwivedi~\cite{JoDw1992}, see also the work of Carr and Gundlach~\cite{CaGu2003}.

The above analysis shows that informally $\mathcal B_1$ is the first outgoing null-geodesic that ``carries information" from the singularity $\mathcal O$ to infinity. In particular, the extension that we have constructed all the way to the massive singularity $\MS$ is in general not unique and other extensions pass $\mathcal B_1$ are in principle possible. Our extension however is unique in the class of exact self-similar extensions - this situation parallels the one in~\cite{Ch1994}.

{\em Step 6. Asymptotic flattening and naked singularities.}

The exactly self-similar RLP spacetime constructed above is, due to simple scaling considerations, not asymptotically flat. In fact, one can easily show that
\[
m(\tau,r) \asymp_{r\to\infty} r,
\]
where $m(\tau,r)$ is the local mass and therefore the total ADM-mass of the RLP-spacetime is infinite. One way to produce an asymptotically flat spacetime
with a naked singularity is to truncate it in the asymptotic region (like in~\cite{Ch1994}), but this is tricky as the Euler-part of the evolution becomes 
highly degenerate in the presence of vacuum. Instead, we resort to a tail dampening procedure so to produce a density profile that decays sufficiently fast at 
null-infinity to give us an asymptotically flat spacetime. 

\begin{figure}

\begin{center}
\begin{tikzpicture}
\begin{scope}[scale=0.5, transform shape]

\coordinate [label=left:$\mathcal O$] (A) at (-2,-2){};

\coordinate (I) at (0,-4);

\coordinate (J) at (4,4);

\coordinate (L) at (-2,-4);

\coordinate (P) at (0,0);

\coordinate (Q) at (2,-2){};

\coordinate (S) at (6,2){};

\node at (-2,-3) {$r=0$};

\node at (0,-2){exact RLP solution};

\node at (0.8,-1.2) {$\uC$};

\node at (4.6,0.3) {$\C$};


\draw [fill=gray!30, dashed] (J) -- (P) -- (Q) -- (S);

\node at (2,0){modified region};

\coordinate [label=left:$\mathcal N$] (F) at (0,-4){};

\draw[thick] (A)--(L);

\draw[thick] (A)--(J);

\draw[thick] (A)--(I);

\draw[thick] (P)--(Q);

\draw[fill=white] (-2,-2) circle (3pt);

\draw[fill=black] (2,-2) circle (2pt);

\draw[fill=black] (0,0) circle (2pt);

\draw[thick] (I)--(S);

\node at (3.8,1.6) {$\D$};

\end{scope}
 \end{tikzpicture}
   \end{center}
      \caption{The self-similar profile is modified inside the grey region $\D$ and the data are given on $\mathcal C$ and $\underline{\mathcal C}$}
     \label{F:CUTOFF}
\end{figure}
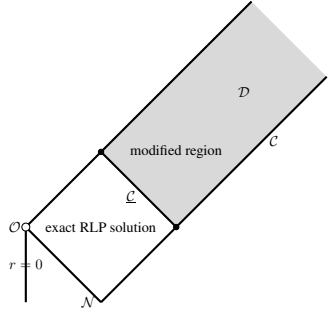

To do this, we employ the so-called double null gauge in which the metric takes the form
\begin{align}\label{E:DNG}
g = - \Omega^2 \, d\u\,d q + r^2 \gamma,
\end{align}
where surfaces of constant $\u$ are outgoing null-surfaces and 
surfaces of constant $q$ are the ingoing null-surfaces. Our goal is to recast Einstein-Euler system in this gauge 
and formulate the corresponding characteristic value problem in a semi-infinite characteristic rectangle depicted by the grey shaded region $\mathcal D$
in Figure~\ref{F:CUTOFF}.
The Einstein field equations in double-null gauge and radial symmetry take the form~\cite{DaRe2005,GHJ2021c} (see also~\cite{BuLe2014} for a formulation in Eddington-Finkelstein coordinates)
\begin{align}
\pa_p \pa_q r & = - \frac{\Om^2}{4r} - \frac{1}{r} \pa_p r \pa_q r  + \pi r\Om^4 T^{p q},  \label{E:DN1}\\
\pa_p \pa_q \log \Om & = - \frac{1+\k}{1-\k}\pi\Om^4 T^{pq} +\frac{\Om^2}{4r^2} + \frac1{r^2} \pa_p r \pa_q r, \label{E:DN2} \\
\pa_q \left(\Om^{-2}\pa_q r\right) & = - \pi r \Om^2 T^{pp}, \label{E:DNCONSTRAINT1}\\
\pa_p\left(\Om^{-2}\pa_p r\right) & = - \pi r \Om^2 T^{qq}.\label{E:DNCONSTRAINT2}
\end{align}
The fluid 4-velocity in the new frame is written as
\[
u = u^p\pa_p + u^q \pa_q
\]
and one can readily compute the components of the energy-momentum tensor:
\begin{align}
& T^{\u\u} = (1+\k)\rho (u^\u)^2, \ \ T^{vv} = (1+\k)\rho (u^q)^2, \ \ T^{\u q} =  
(1-\k)\rho \Om^{-2},  \label{E:TPP}\\ 
& T^{AB} 
= \k \rho r^{-2} \gamma^{AB},   \ \
 T^{\u A} = T^{ q A} = 0, \ \ A,B = 2,3. \label{E:TAB} 
\end{align}
The normalisation condition~\eqref{E:NORMALISATION} now reads
\begin{align}
T^{\u\u} T^{ q q}
= \left(\frac{1+\k}{1-\k}\right)^2(T^{\u q})^2 . \label{E:TCONSTRAINT}
\end{align}
To complete the formulation of the problem, we must explain how the fluid unknowns $(\rho,u^p,u^q)$ evolve in time. 
The corresponding Euler evolution is a simple consequence of the Bianchi identities and reads:
\begin{align}
\pa_p (\Om^4r^2T^{pp}) + \frac{\Om^2}{r^{\frac{4\k}{1-\k}}} \pa_q ( \Om^2r^{\frac{2+2\k}{1-\k}} T^{pq}  ) &=0, \label{E:DNEULER1}\\
\pa_p ( \Om^2r^{\frac{2+2\k}{1-\k}} T^{pq}  )  + \frac{r^{\frac{4\k}{1-\k}}}{\Om^2}\pa_q (\Om^4r^2T^{qq}) & = 0.\label{E:DNEULER2}
\end{align}
Equations~\eqref{E:DNEULER1}--\eqref{E:DNEULER2} are nonlinear transport equations and the full coupling~\eqref{E:DN1}--\eqref{E:DNEULER2}
is a nonlinear system of coupled wave- and transport-like equations. By the finite speed of propagation, we may hope to modify the RLP-solution in an asymptotic region
in such a way that the solution is exactly self-similar and agrees with the constructed RLP-solution in a region in the causal past of some ingoing curve $\uC$ (the white region in Figure~\ref{F:CUTOFF}), and otherwise solve a characteristic initial value problem in the gray shaded region $\mathcal D$. To make this work, we prescribe the data along the two characteristic characteristic curves $\C$ and $\uC$, so that the constraint equations~\eqref{E:DNCONSTRAINT1}--\eqref{E:DNCONSTRAINT2} are satisfied respectively.  

The final and key technical step is the proof of local well-posedness for this characteristic value problem. Although various local well-posedness results for the Einstein-Euler system are availalble (see for example~\cite{CB,Rendall1992, EhKi1993,BrKa} and references therein), we have not been able to find such a reference that treats isothermal EE-system in the double-null coordinate gauge. To do this, we make a change of fluid unknowns that effectively ``diagonalises"~\eqref{E:DNEULER1}--\eqref{E:DNEULER2} in the spirit of Riemann variables in the theory of hyperbolic conservation laws. We do not present the details (an interested reader can see Section 9.1 of~\cite{GHJ2021c}), but the key geometric feature of this reformulation is to highlight the relationship of the two characteristic geometries in the problem - the light and the acoustical null-cones, see Figure 9 in~\cite{GHJ2021c}.

Finally, we impose data on the outgoing piece of the boundary $\C$ that are asymptotically flat and belong to a suitable weighted Sobolev space, while the data on the ingoing piece coincide precisely with the self-similar RLP-solution. Using the method of characteristics we prove local-in-$\delta$ well-posedness of the characteristic value problem, where $\delta$ is the width of the semi-infinite characteristic rectangle $\mathcal D$. The resulting spacetime, obtained by gluing the solution inside $\mathcal D$ with the exact self-similar solution is a classical solution by the finite speed of propagation with a Penrose diagram given by~Figure~\ref{F:THMPIC2}. It is then not hard to show that the Cauchy horizon of the resulting spacetime is but a suitable ``deformation" of $\mathcal B_1$ from Figure~\ref{F:THMPIC1} and that it contains a naked singularity in the sense of Definition 1.1 in~\cite{RoSR2019}, see also Figure~\ref{F:NS}.

\end{document}